
  %
  %
  %



  \newcount\fontset
  \fontset=1
  \def\dualfont#1#2#3{\font#1=\ifnum\fontset=1 #2\else#3\fi}
  \dualfont\eightrm {cmr8} {cmr7}
  \dualfont\eightsl {cmsl8} {cmr7}
  \dualfont\eightit {cmti8} {cmti10}
  \dualfont\eightmi {cmmi8} {cmmi7}
  \dualfont\eightbf {cmbx8} {cmbx10}
  \dualfont\fivemi {cmmi5} {cmmi7}
  \dualfont\tensc {cmcsc10} {cmcsc10}
  \dualfont\eighttt {cmtt8} {cmtt10}
  \dualfont\tenss {cmss10} {cmr10}
  \dualfont\titlefont {cmbx12} {cmbx10}
  \dualfont\eightsymbol {cmsy8} {cmsy10}


  \def\vg#1{\ifx#1\null\null\else
    \ifx#1\ { }\else
    \ifx#1,,\else
    \ifx#1..\else
    \ifx#1;;\else
    \ifx#1::\else
    \ifx#1''\else
    \ifx#1--\else
    \ifx#1))\else
    { }#1\fi \fi \fi \fi \fi \fi \fi \fi \fi}
  \newcount\secno \secno=0
  \newcount\stno
  \def\goodbreak{\vskip0pt plus.1\vsize\penalty-250
    \vskip0pt plus-.1\vsize\bigskip}
  \outer\def\section#1{\stno=0
    \global\advance\secno by 1
    \goodbreak\vskip\parskip
    \message{\number\secno.#1}
    \noindent{\bf\number\secno.\enspace#1.\enspace}}
  \def\seqnumbering{\global\advance\stno by 1
    \number\secno.\number\stno}
  \def\label#1{\global\edef#1{\number\secno.\number\stno}}
  \def\sysstate#1#2#3{\medbreak\noindent {\bf\seqnumbering.\enspace
    #1.\enspace}{#2 #3\vskip 0pt}\medbreak}
  \def\state#1 #2\par{\sysstate{#1}{\sl}{#2}}
  \def\definition#1\par{\sysstate{Definition}{\rm}{#1}}
  \def\proof{\medbreak\noindent{\it Proof.\enspace}}
  \def\proofend{\ifmmode\eqno\square\else\hfill\square\medbreak\fi}
  \newcount\zitemno \zitemno=0
  \def\zitem{\global\advance\zitemno by 1 \smallskip
    \item{\ifcase\zitemno\or i\or ii\or iii\or iv\or v\or vi\or vii\or
    viii\or ix\or x\or xi\or xii\fi)}}
  \def\$#1{#1 $$$$ #1}


  \def\({\left(}
  \def\){\right)}
  \def\[{\left\Vert}
  \def\]{\right\Vert}
  \def\*{\otimes}
  \def\+{\oplus}
  \def\:{\colon}
  \def\<{\langle}
  \def\>{\rangle}
  \def\and{\hbox{\quad and \quad}}
  \def\arw{\rightarrow}
  \def\calcat#1{\,{\vrule height8pt depth7pt}_{\,#1}}
  
  \def\crossproduct{\hbox to 1.8ex{$\times \kern-.45ex\vrule height1.1ex
depth0pt width0.45truept$\hfill}}
  \def\cstar{$C^*$}
  \def\for#1{,\quad \hbox{ for }#1} 
  \def\inv{^{-1}}
  \def\pmatrix#1{\left[\matrix{#1}\right]}
  
  \def\square{\hbox{$\sqcap\!\!\!\!\sqcup$}}
  \def\stress#1{{\it #1}\/}
  
  \def\x{\times}
  \def\|{\Vert}
  \def\inv{^{-1}}


  \newcount\bibno \bibno=0
  \def\newbib#1{\global\advance\bibno by 1 \edef#1{\number\bibno}}
  \def\cite#1{{\rm[\bf #1\rm]}}
  \def\scite#1#2{[{\bf #1},{\it #2\/}]}
  \def\lcite#1{#1}
  \def\se#1#2#3#4{\def\a{#1}\def\b{#2}\ifx\a\b#3\else#4\fi}
  \def\setem#1{\se{#1}{}{}{, #1}}
  \def\index#1{\smallskip \item{[#1]}}
  \def\tit#1{``#1''}

  \def\zarticle#1, auth = #2,
   title = #3,
   journal = #4,
   year = #5,
   volume = #6,
   pages = #7,
   NULL#8 xyzzy {\index{#1} #2, \tit{#3}, {\sl #4\/} {\bf #6} (#5), #7.}

  \def\ztechreport#1,
    auth = #2,
    title = #3,
    institution = #4,
    year = #5,
    type = #6,
    note = #7,
    NULL#8
    xyzzy
    {\index{#1} #2, \tit{#3}\setem{#6}\setem{#4}\setem{#5}\setem{#7}.}

  \def\zunpublished#1,
    auth = #2,
    title = #3,
    institution = #4,
    year = #5,
    type = #6,
    note = #7,
    NULL#8
    xyzzy
    {\index{#1} #2, \tit{#3}\setem{#6}\setem{#4}\setem{#5}\setem{#7}.}

  \def\zbook#1,
    auth = #2,
    title = #3,
    publisher = #4,
    year = #5,
    volume = #6,
    series = #7,
    NULL#8
    xyzzy
    {\index{#1} #2, \tit{#3}\setem{#7}\se{#6}{}{}{ vol. #6}, #4, #5.}

  \def\zmasterthesis#1,
    auth = #2,
    title = #3,
    school = #4,
    year = #5,
    type = #6,
    NULL#7
    xyzzy
    {\index{#1} #2, \tit{#3}, #6, #4, #5.}

  \def\zinproceedings#1,
  auth = #2,
    title = #3,
    booktitle = #4,
    year = #5,
    pages = #6,
    organization = #7,
    note = #8,
    NULL#9
    xyzzy
    {\index{#1} #2, \tit{#3}\se{#4}{}{}{, In
      {\sl #4}}\setem{#5}\setem{#6}\setem{#7}\setem{#8}.}

  \def\zphdthesis#1,
    auth = #2,
    title = #3,
    school = #4,
    year = #5,
    type = #6,
    NULL#7
    xyzzy
    {\index{#1} #2, \tit{#3}, #6, #4, #5.}

  \def\zbooklet#1,
    auth = #2,
    title = #3,
    howpublished = #4,
    year = #5,
    NULL#6
    xyzzy
    {\index{#1} #2, \tit{#3}, #4, #5.}

  \def\zmisc#1,
    auth = #2,
    title = #3,
    note = #4,
    howpublished = #5,
    year = #6,
    NULL#7
    xyzzy
    {\index{#1} #2, \tit{#3}\setem{#4}\setem{#5}\setem{#6}.}


  \nopagenumbers
  \voffset=2\baselineskip
  \advance\vsize by -\voffset
  \headline{\ifnum\pageno=1 \hfil \else \tensc\hfil
    deformation quantization via fell bundles
  \hfil\folio \fi}


  \def\pmatrix#1{\left[\matrix{#1}\right]}
  \def\={:=}
  \def\vaipara{\ \longmapsto\ }
  \def\DownParenthesisfill{$
    \braceld\leaders\vrule\hfill\bracerd$}
  \def\OverParenthesis#1{\mathop{\vbox{\ialign{##\crcr\noalign{\kern3pt}
    \DownParenthesisfill\crcr\noalign{\kern3pt\nointerlineskip}
    $\hfil\displaystyle{#1}\hfil$\crcr}}}\limits}
  \def\class#1{\OverParenthesis{#1}^{.}}
  \def\pdop#1#2{\partial_{#1}\!(#2)}
  \def\B{{\cal B}}
  \def\D{{\cal D}}
  \def\C{{\bf C}}
  \def\Z{{\bf Z}}
  \def\T{{\bf T}}
  \def\R{{\bf R}}
  \def\defprod{\times}
  \def\defstarsymbol{\diamond}
  \def\defstar{^\defstarsymbol}
  \def\inv{^{-1}}
  \def\t{\theta}
  \def\g{\gamma}
  \def\deriv#1#2#3{{d \over d#2}\bigl(#1\bigr)\calcat{#2=#3}}
  \def\action{\phi}
  \def\DB#1#2{{#1}^{#2}}
  \def\DA#1#2#3{{#1}^{#3}_{#2}}
  \def\DBs{\DB{\B}{\t}}
  \def\DAs{\DA{B}{\g}{\t}}
  \def\Gh{\widehat G}
  \def\dd{deformation data\vg}
  \def\dds{deformation data\vg}
  \def\tt{\vartheta}
  \def\expi#1{{\rm e}^{2\pi i #1}}
  \def\expineg#1{{\rm e}^{-2\pi i #1}}


  \newbib\AEE
  \newbib\Picard
  \newbib\Brat
  \newbib\Dixmier
  \newbib\Soft
  \newbib\Circle
  \newbib\Amena
  \newbib\FD
  \newbib\QCan
  \newbib\Mt
  \newbib\MatsumotoHopf
  \newbib\MaTo
  \newbib\Ng
  \newbib\Irrat
  \newbib\RieffelCF
  \newbib\RieffelHeis
  \newbib\RieffelDQi
  \newbib\RieffelDQii
  \newbib\Woro
  \newbib\ZM


  \null
  \vskip-2\bigskipamount
  \centerline{\titlefont DEFORMATION QUANTIZATION VIA}
  \smallskip
  \centerline{\titlefont FELL BUNDLES}
  \bigskip
  \centerline{\tensc
  Beatriz Abadie\footnote{*}{\eightrm Partially supported by CONICYT,
Proyecto 2002 -- Uruguay.},
  and
  Ruy Exel\footnote{**}{\eightrm Partially supported by CNPq --
Brazil.}}
  \bigskip \medskip
  \centerline{\bf May 1, 1997}


  \bigskip
  \bigskip
  \midinsert
  \narrower
  \narrower
  \baselineskip=9.5pt
  \eightbf
  \noindent
  ABSTRACT. A method for deforming C*-algebras is introduced, which
applies to C*-algebras that can be described as the cross-sectional
C*-algebra of a Fell bundle.  Several well known examples of
non-commuta\-ti\-ve algebras, usually obtained by deforming commutative
ones by various methods, are shown to fit our unified perspective of
deformation via Fell bundles.  Examples are
  the non-commutative spheres of Matsumoto,
  the non-commutative lens spaces of Matsumoto and Tomiyama,
  and
  the quantum Heisenberg manifolds of Rieffel.
  In a special case, in which the deformation arises as a result of an
action of R$^{\fivemi 2d}$, assumed to be periodic in the first d
variables, we show that we get a strict deformation quantization.
  \endinsert

  \section{Introduction}
  Deformations of \cstar-algebras, specially of the \cstar-algebra of
continuous functions on the phase space of a classical physical system,
have long been associated to the process of quantization and have been
used to explain quantum mechanical phenomena such as the correspondence
principle (see, for example, \cite{\RieffelHeis}).

One of the most popular processes for constructing these deformations is
to describe the given \cstar-algebra by means of generators and
relations, and, after introducing a deformation parameter into these
relations, to consider the universal \cstar-algebra for the new
relations.  This process can be used, for example, for constructing
  the non-commutative torus \cite{\Irrat},
  the soft torus \cite{\Soft},
  the quantum $SU_2$ groups \cite{\Woro},
  the non-commutative spheres \cite{\Mt},
  the non-commutative lens spaces \cite{\MaTo},
  and
  the algebra of the $q$-canonical commutation relations \cite{\QCan}.

However, \cstar-algebras arising from generators and relations are often
intractable objects motivating one to search for alternate
constructions.  The goal of the present work is to show the usefulness
of the techniques of Fell bundles (also known as \cstar-algebraic
bundles \cite{\FD}) in the study of deformations of \cstar-algebras.

The first step in our construction, as our title suggests, requires one
to look for a Fell bundle structure for the given algebra $B$, over a
locally compact topological group $G$.  In the important special case of
a discrete group $G$, by far the easiest to understand, this is roughly
equivalent
  (see \cite{\Amena})
  to finding a $G$-grading for $B$, that is, a decomposition of $B$ as
the closure of the direct sum $\bigoplus_{t \in G} B_t$ of a family of
closed linear subspaces
  $\B = \{B_t\}_{t\in G}$,
  satisfying
  $B_t B_s \subseteq B_{ts}$, and
  $B_t^* = B_{t\inv}$ for all $t$ and $s$ in $G$.
  The collection
  $\B = \{B_t\}_{t\in G}$, equipped with the operations of
multiplication
  $$
  \cdot : B_t \x B_s \arw B_{ts},
  $$
  and involution
  $$
  * : B_t \arw B_{t\inv},
  $$
  inherited from $B$, gives an important example of a Fell bundle.  This
Fell bundle contains, in many cases, the necessary information for the
reconstruction of the algebra $B$.  This is so, for example, when the
group $G$ is amenable.  However, in some important cases where $G$ does
not posses this property, this is still true.  See \cite{\Amena} and
\cite{\Ng} for a thorough discussion of this problem.

The second step in constructing our deformation requires an action $\t$
of the group $G$ on $\B$, called the \stress{deforming} action, which is
then used to deform the Fell bundle structure by means of introducing a
new multiplication operation $\defprod$ and a new involution $\defstar$,
via the formulas
  $$
  a_t\defprod b_s \= a_t \t_t(b_s),
  $$
  and
  $$
  a_t\defstar \= \t_t\inv(a_t^*),
  $$
  for $s,t\in G$, $a_t$ in $B_t$ and $b_s$ in $B_s$.

The norm and the linear structure of $\B$, on the other hand, are kept
intact.  The deformed algebra is then obtained by taking the
cross-sectional \cstar-algebra of the deformed Fell bundle.

The invariance of the linear structure and the norm of the Fell bundle,
under the deformation, is an important feature of our construction,
because it allows us to embed part of the original algebra into its
deformed version.

An important ingredient in the established theories of deformations of
\cstar-algebras is the continuity of the collection of deformed
algebras, as a function of the deformation parameter.  This property is
usually expressed by the fact that these form a continuous field of
\cstar-algebras over the parameter space (see \cite{\RieffelCF}).  Our
method of deformation is also shown to produce continuous fields of
\cstar-algebras when one is given, not just one, but a family
$\{\theta^{\hbar}\}_{\hbar\in I}$ of actions of $G$ on $\B$, as above,
where $I$ is an interval of real numbers.  Since all of the applications
presented here deal with discrete groups, we have opted to restrict our
study of continuity for the discrete group case, although it seems
plausible to expect that generalizations can be found for continuous
groups.

  The continuity results we obtain are, essentially, reworkings of
Rieffel's ideas in \cite{\RieffelCF} for our more general situation of
Fell bundles.

Our process of deformation is rather simple to follow, given that the
necessary ingredients, the Fell bundle structure and the deforming
action $\t$, are provided.  One way to obtain these ingredients is via a
\stress{\dd{}} for $B$, that is, a triple
  $(G,\g,\t)$,
  where $G$ is a discrete abelian group and $\g$ and $\t$ are commuting
actions on $B$, respectively of the Pontryagin dual $\Gh$ of $G$, and of
$G$ itself.  Since $\Gh$ is compact, the spectral decomposition of $\g$
provides the Fell bundle structure, while $\t$ plays the role of the
deforming action.
  This approach essentially consists of introducing a deformation
parameter after taking a certain Fourier transform, a method that has
already been used by other authors, including Rieffel (see, for example,
the formula for the definition of $*_\hbar$ on page 541 of
\cite{\RieffelHeis}).  The advantage of emphasizing the Fell bundle
structure is, perhaps, in making some formulas more transparent.

  This construction, albeit rather elementary, provides some very
interesting examples.  We show, for instance, that
  the non-commutative spheres,
  the non-commutative lens spaces,
  and
  the quantum Heisenberg manifolds \cite{\RieffelHeis}, can all be seen
under the unified perspective of deformation of Fell bundles.

  Even though our deformation is described for Fell bundles over general
locally compact groups, all of our examples are restricted to the
simpler case of discrete groups.  This is partly due to the difficulty
in identifying Fell bundle structures over non-discrete groups, but we
feel that more general examples may be just as relevant.

We have already touched upon the need to search for alternate methods of
deformations of \cstar-algebras.  Perhaps one of the most successful
approaches is the one due to Rieffel \cite{\RieffelDQi,\RieffelDQii}
where deformations arise from actions of $\R^n$ on the given algebra.
As indicated by Rieffel in page 84 of \cite{\RieffelDQii}, when the
action of $\R^n$ factors through a compact group, one can describe the
algebra as the cross-sectional algebra of the Fell bundle arising from
the spectral subspaces of the given action.  Rieffel's deformation may
then be seen as a deformation of the Fell bundle structure by means of
the introduction of a 2-cocycle.

  Roughly following this approach, our method may be viewed as a study
of the situations in which one benefits from the simplifications arising
{}from the general phenomena of ``compactness'', or the dual notion of
``discreteness''.
  However, our approach differs from Rieffel's because our deformation
is caused by a group action, as opposed to a 2-cocycle, and because we
deform both the multiplication operation and the involution, while
Rieffel's deformation affects only the former.

  One of the most important ingredients in the theory of deformation
quantization is the computation of the first order term (the derivative
at zero) of the deformed multiplication operation, as a function of the
deformation parameter, often denoted by $\hbar$.  We carry this out in
the case where a \cstar-algebra $B$ is deformed via a \dd arising from
an action $\action$ of $\R^{2d}$ on $B$.  To be precise, let $\action$
be such an action.  Our ``compactness'' assumption consists of supposing
that $\action$ is periodic in its first $d$ variables, hence inducing an
action $\g$ of the compact $d$-dimensional torus $\T^d$ on $B$ .  Fixing
a real number $\hbar$, we let $\t^\hbar$ be the action of $\Z^d$ on $B$
given by
  $$
  \t^\hbar_{(n_1,\ldots,n_d)} =
  \action_{(0,\ldots,0,\hbar n_1,\ldots,\hbar n_d)}
  \for (n_1,\ldots,n_d) \in \Z^d.
  $$
  The triple $(\Z^d,\g,\t^\hbar)$ is then a \dd for $B$, with which we
construct the deformed algebra $B^{(\hbar)}$.  We then show that
  $$
  \lim_{\hbar \rightarrow 0}
  \left\Vert {f \defprod_\hbar g - fg \over \hbar} -
  {1 \over 2\pi i}
  \sum_{j=1}^d \pdop{x_j}{f} \pdop{y_j}{g} \right\Vert_\hbar
  = 0,
  $$
  where $\defprod_\hbar$ and $\|\cdot\|_\hbar$ refer to the deformed
product and norm of $B^{(\hbar)}$, respectively, and $\partial_{x_j}$
and $\partial_{y_j}$ denote the derivation operators associated to the
action $\action$ under the coordinate system
$(x_1,\ldots,x_d,y_1,\ldots,y_d)$ for $\R^{2d}$.

This is initially done for a very restrictive class of elements $f$ and
$g$ in $B$, namely the smooth elements belonging, each, to a single
spectral subspace for $\g$.  The proof of this result is extremely
simple and the formulas involved show, in a very transparent way, the
roles of the various ingredients present in the context.  In particular,
the heavy machinery of oscillatory integrals of \cite{\RieffelDQi} does
not intervene, thanks, of course, to the simplification introduced by
the periodicity assumption.  Because of the simplicity of our formulas,
our approach may be pedagogically relevant for the understanding of more
sophisticated constructions.

The formula for the derivative of the deformed multiplication, above, is
then extended to smooth elements $f$ and $g$, with a proof which permits
one to see through the extent to which the differentiability properties
of these elements are necessary.

This immediately implies that
  $$
  \lim_{\hbar\rightarrow 0}
  \left\Vert {f \defprod_\hbar g - g \defprod_\hbar f - [f,g]
  \over \hbar} -
  {1 \over 2\pi i}
  \{f,g\}
  \right\Vert_\hbar = 0,
  $$
  where $\{\cdot,\cdot\}$ is the Poisson bracket induced by $\action$,
and
  $[\cdot,\cdot]$ is the commutator for the original multiplication on
$B$.

  Combining this with the fact, shown below, that the $B^{(\hbar)}$,
form a continuous field of \cstar-algebras, we get a strict deformation
quantization in the sense of Rieffel \cite{\RieffelDQi,\RieffelDQii}.

  The authors would like to acknowledge the support of CONICYT (Uruguay)
and FAPESP (Brazil) for funding numerous academic visits while this
research was conducted.

  \section {The Deformation}
  Let $G$ be a locally compact topological group and let
  $\B = \{B_t\}_{t\in G}$ be a
  \cstar-algebraic bundle over $G$.  The reader is referred to
\cite{\FD} for a comprehensive treatment of the basic theory of
  \cstar-algebraic bundles.
  These objects have recently been referred to as ``Fell bundles'', a
terminology we have chosen to adopt.  In what follows, we shall identify
$\B$ with the total bundle space
  $
  \bigcup_{t\in G} B_t
  $.

  Let $\D = \{D_t\}_{t\in G}$ be another Fell bundle over $G$.  A map
$\psi$ from $\B$ to $\D$
  is called a \stress{homomorphism} if
  \zitemno = 0
  \zitem $\psi$ is continuous,
  \zitem $\psi(B_t) \subseteq D_t$, for all $t$ in $G$,
  \zitem $\psi$ is linear on each $B_t$,
  \zitem $\psi(a b) = \psi(a) \psi(b)$, for all $a,b$ in $\B$, and
  \zitem $\psi(a^*) = \psi(a)^*$, for all $a$ in $\B$.

  \medskip
  Let $\psi$ be a homomorphism from $\B$ to $\D$.  Observe that, since
$\psi$ restricts to a *-homomorphism between the \cstar-algebras $B_e$
and $D_e$
  (where $e$ denotes the unit group element), then it is necessarily
contractive there.  Also, for each $b_t$ in $B_t$ we have
  $$
  \|\psi(b_t)\|^2 =
  \|\psi(b_t^*b_t)\| \leq
  \|b_t^*b_t\| =
  \|b_t\|^2,
  $$
  so that $\psi$ is in fact norm-contractive everywhere.

  If $\psi$ is bijective, then $\psi\inv$ is continuous as well
  \scite{\FD}{II.13.17}
  and hence it is also a homomorphism.  In this case we say that $\psi$
is an \stress{isomorphism}.  If, in addition, $\D = \B$, then $\psi$ is
called an \stress{automorphism} of $\B$.
  In particular, if $\psi$ is an isomorphism, then it must be isometric.
  See also \scite{\FD}{VIII.3.3}.

  Given another locally compact topological group $H$, by an
\stress{action} of $H$ on $\B$, we shall mean an assignment $\t$ which,
for each $x$ in $H$, gives an automorphism $\t_x$ of $\B$, satisfying
  $\t_x \t_y = \t_{xy}$, for $x,y \in H$.
  We shall say that $\t$ is a continuous action if the map
  $$
  (x,b) \in H \x \B \vaipara \t_x(b) \in \B
  $$
  is jointly continuous.

  Let us now suppose we are given a Fell bundle $\B = \{B_t\}_{t\in G}$
over the locally compact group $G$, as well as a continuous action $\t$
of the very same group $G$ on $\B$.  We wish to construct a new product
on $\B$, denoted $\defprod$, and a new involution, called $\defstar$,
providing a ``deformed'' bundle structure.  In order to do so, define
for
  $a_t$ in $B_t$ and $b_s$ in $B_s$,
  $$
  a_t\defprod b_s = a_t \t_t(b_s),
  $$
  and
  $$
  a_t\defstar = \t_t\inv(a_t^*).
  $$

  \state Proposition If $\B$ keeps its linear, topological and norm
structure, but is given the deformed operations $\defprod$ and
$\defstar$, then it is a Fell bundle.

  \proof To check that the new multiplication operation is continuous,
we shall use \scite{\FD}{VIII.2.4}.  That is, given continuous sections
$\beta$ and $\g$ of $\B$, we must show that the map
  $$
  (r,s)\in G \x G \vaipara \beta(r)\defprod\g(s)\in \B
  $$
  is continuous.  Now, we have
  $
  \beta(r)\defprod\g(s) = \beta(r) \t_r(\g(s)),
  $
  which is continuous by the continuity of $\t$ and of the original
multiplication.  A similar argument shows that the deformed involution
is continuous.

  Let us now verify the associativity of $\defprod$.  Given $a_r$ in
$B_r$, $b_s$ in $B_s$, and $c_t$ in $B_t$ we have
  $$
  (a_r \defprod b_s) \defprod c_t =
  (a_r \t_r(b_s) ) \defprod c_t =
  a_r\t_r(b_s)\t_{rs}(c_t)
  \$=
  a_r\t_r(b_s\t_s(c_t)) =
  a_r \defprod ( b_s \t_s(c_t)) =
  a_r \defprod ( b_s \defprod c_t).
  $$

  As for the anti-multiplicativity of the involution, let $a_r\in B_r$
and $b_s\in B_s$.  Then
  $$
  (a_r \defprod b_s)\defstar
  =
  (a_r \t_r(b_s) )\defstar
  =
  \t_{(rs)\inv} (a_r \t_r(b_s) )^*
  =
  \t_{s\inv} \t_{r\inv} (\t_r(b_s^*) a_r^* )
  \$=
  \t_{s\inv}(b_s^*) \t_{s\inv} \t_{r\inv} (a_r^*)
  =
  \t_{s\inv}(b_s^*) \defprod \t_{r\inv} (a_r^*)
  =
  b_s\defstar \defprod a_r\defstar.
  $$

  The verification of the remaining axioms is routine and so is left as
an exercise.
  \proofend

  \definition The bundle constructed above will be called the
$\t$-deformation of $\B$ and will be denoted $\DBs$.

  Recall that a Fell bundle is said to be \stress{saturated}
\scite{\FD}{VIII.2.8} if $B_{rs} = B_r B_s$ (closed linear span) for all
$r,s$.  In the special case that $G$ is equipped with a ``length''
function
  $$
  |\cdot| : G \arw \R_+
  $$
  satisfying $|e| = 0$, and the triangular inequality
  $|rs| \leq |r|+|s|$,
  then we say that $\B$ is \stress{semi-saturated}
  (see \scite{\Circle}{4.1, 4.8}, \scite{\Amena}{6.2}),
  if $B_{rs} = B_r B_s$, whenever $r,s \in G$ are such that
  $|rs| = |r|+|s|$.

  \state Proposition \label{\Satura} If $\B$ is saturated (resp.
semi-saturated) then so is $\DBs$.

  \proof It is enough to observe that
  $B_r \defprod B_s = B_r \t_r(B_s) = B_r B_s$.
  \proofend

  \begingroup
  \def\CE{{\cal E}}
  \def\stackrel#1#2{\buildrel#1 \over #2}
  \def\frac#1#2{{#1 \over #2}}
  \def\lrarw{\rightarrow}
  \def\ccb{\bigoplus_{x\in G}B_x}
  \def\RR{{\bf R}}
  
  \def\prodh{\defprod_{\hbar}}
  \def\starh{^{\defstarsymbol_\hbar}}
  \def\tili{\tilde{i}}
  \section {Continuous fields arising from deformations}
  \global\edef\ContSection{\number\secno}
  The purpose of this section is to show that the collection of deformed
algebras, originated from a continuous family of group actions on a Fell
bundle, gives rise to a continuous field of \cstar-algebras.

We first establish some facts on Fell bundles over discrete groups that
will enable us to extend the techniques in \cite{\RieffelCF} to discuss
upper semicontinuity.  Let $\B$ and $\D$ be fell bundles over a discrete
group $G$, and let $\Phi: {\cal{D}}\lrarw {\cal{B}}$ be a Fell bundle
homomorphism. Since $\Phi$ is contractive, one can define $\Phi^1:
L^1({\cal{D}})\lrarw L^1({\cal{B}})$ by $[\Phi(f)](x)=\Phi[f(x)]$, for
$f\in L^1({\cal{D}})$, and $x\in G$. It is easily checked that $\Phi^1$
is a $*$-algebra homomorphism, so it gives rise to a \cstar-algebra
homomorphism ${\tilde{\Phi}}: C^*({\cal{D}})\lrarw C^*({\cal{B}})$.

A sequence of Fell bundle homomorphisms
  $$
  0\lrarw {\cal{E}}
\stackrel{i}{\lrarw}{\cal{D}}\stackrel{\Pi}{\lrarw}{\cal{B}}\lrarw 0
  $$
  is said to be exact if so are the sequences
  $$
  0\lrarw \CE_x \stackrel{i|_{\CE_x}}{\lrarw}
\D_x\stackrel{\Pi|_{\D_x}}{\lrarw}\B_x\lrarw 0
  $$
  for all $x\in G$.

  \state Lemma \label{\exs} Let
  $ 0\lrarw {\cal{E}} \stackrel{i}{\lrarw}{\cal{D}}
  \stackrel{\Pi}{\lrarw}{\cal{B}}\lrarw 0$
  be an exact sequence of Fell bundle homomorphisms over a discrete
group $G$.
  Then $ 0\lrarw C^*(\CE) \stackrel{\tili}{\lrarw}
C^*(\D)\stackrel{\tilde{\Pi}}{\lrarw}C^*(\B)\lrarw 0$ is also exact.

  \proof In view of \scite{\ZM}{2.29}, and \scite{\FD}{VIII 5.11, 16.3},
we only need to show that $ 0\lrarw L^1(\CE) \stackrel{i^1}{\lrarw}
L^1(\D)\stackrel{\Pi^1}{\lrarw} L^1(\B)\lrarw 0$ is exact.  It is
apparent from the definition that $i^1$ is injective, and that
Im$(i^1)=$ker$(\Pi^1)$, so we need only show that $\Pi^1$ is onto.
 Fix $b_x\in \B_x$ and $\epsilon >0$. Since $b_x\delta_x\in$
Im$\tilde{\Pi}$, there exists $\tilde{d}\in C^*(\D)$ such that
$\tilde{\Pi}(\tilde{d})=b_x\delta_x$, and
  $$
  \|\tilde{d}\|_{C^*(\D)}\leq \|\tilde{d}+
\hbox{ker}\tilde{\Pi}\|_{C^*(\D)/\hbox{ker}\tilde{\Pi}}+\epsilon
=\|b_x\delta_x\|_{C^*(\B)}+\epsilon=\|b_x\|_{\B_x}+\epsilon.
  $$
  
Let $P_x^{\D}$ (resp. $P_x^{\B}$) denote the projection onto the
$x^{th}$ spectral subspace of $\D$ (resp. $\B$). Then
$P_x^{\B}\tilde{\Pi}=\tilde{\Pi}P_x^{\D}$, since the equality holds when
restricted to $L^1(\D)$. Now set $d= P^{\D}_x(\tilde{d})$.  Then $d\in
\D_x$, $\tilde
\Pi(d)=\tilde{\Pi}P_x^{\D}(\tilde{d})=P^{\B}_x\tilde{\Pi}(\tilde{d})=b_x$,
and $\|d\|_{\D_x}\leq \|\tilde{d}\|_{C^*(\D)}\leq
\|b_x\|_{\B_x}+\epsilon$.

Now, if $\sum b_n\delta_{x_{n}}\in L^1(\B)$, choose as above, for each
positive integer $n$, $d_n\in \D_{x_{n}}$ so that $\Pi(d_n)=b_n$, and
$\|d_n\|_{\D_{x_{n}}}\leq \|b_n\|_{\B_{x_{n}}} +n^{-2}$. Then
$\Pi^1(\sum c_n \delta_{x_{n}})= \sum b_n \delta_{x_{n}}$. So $\Pi^1$ is
onto.  \proofend

Back to the setting of the previous section, we consider a
\cstar-algebra $B$ that can be viewed as the cross-sectional
\cstar-algebra of a Fell bundle $\B$ over a discrete group $G$ whose
$x^{th}$ fiber we denote by $B_x$. At this point we are ready to get a
deformed version of $B$ by means of an action $\theta$ of $G$.

  Notice that the algebra $B^{\theta}$ contains as a dense
$*$-subalgebra the set $\ccb$ of compactly supported
cross-sections. Although the $*$-algebra structure of $\ccb$ depends on
$\theta$, its vector space structure does not.

   Our purpose is to produce a continuous field of \cstar-algebras,
given a family $\{\theta^{\hbar}\}$ of actions of $G$ on $\B$.  The
crucial point is to show that the map $\hbar\mapsto \|\phi\|_{\hbar}$ is
continuous for any $\phi\in \ccb$, where $\|\phi\|_{\hbar}$ denotes the
norm of $\phi$ as an element of $C^*(\B^{\theta^{\hbar}}).$

  \sysstate{Notation}{\rm}
  {\label{\defi} In the context above, let $I\subset \RR$ be an open
interval containing $0$ and, for each $\hbar\in I$, let $\theta^{\hbar}$
be an action of $G$ on the Fell bundle $\B$ such that $\theta^0$ is the
identity, and that the map $\hbar\mapsto \theta^{\hbar}_x(b)$ is
continuous for any fixed $x\in G$, $b\in \B$.  We denote the bundle
$\B^{\theta^{\hbar}}$ by $\B^\hbar$, and by $\prodh$, $\starh$ its
product and involution, respectively. The norm in $C^*(\B^{\hbar})$ is
denoted by $\|$ $\|_{\hbar}$.
 }

  \state Proposition The map $\hbar\mapsto \|\phi\|_{\hbar}$ is upper
semicontinuous on $I$ for all $\phi\in \ccb$.
  
  \proof The proof follows the lines of \cite{\RieffelCF}. Let $\D$ be
the Fell bundle over $G$ whose $x^{th}$ fiber is the Banach space
$D_x=C_0(I,B_x)$, with multiplication and involution given by
  $$
  (f_x\star f_y)(\hbar) =
  f_x(\hbar)\prodh f_y(\hbar)\hbox{,\ \ \ \ }
  f_x^{\star}(\hbar)=(f_x(\hbar))\starh,
  $$
  for $f_x\in D_x$, $f_y\in D_y$. For each $\hbar\in I$ consider the
Fell bundle homomorphism
  $\Pi^{\hbar}: \D \lrarw \B$,
  given by $\Pi^{\hbar}(f)=f(\hbar)$.  Since $\Pi^{\hbar}$ is onto for
any $\hbar\in I$ we get, as in Lemma \lcite{\exs}, the exact sequence
  $$
  0\lrarw C^*(\CE^{\hbar})\stackrel{\tili^{\hbar}}{\lrarw}
C^*(\D)\stackrel{\tilde{\Pi^{\hbar}}}{\lrarw} C^*(\B^{\hbar})\lrarw 0,
  $$
  where $\CE^{\hbar}$ is the Fell bundle whose $x^{th}$ fiber is
$E^{\hbar}_x=$ ker$\Pi_x^{\hbar}$, with the structure inherited from
$\D$, and $\tili^{\hbar}$ denotes inclusion.

In order to apply \scite{\RieffelCF}{1.2}, we next consider $C_0(I)$ as
a \cstar-subalgebra of the algebra of multipliers of $D_e$, in the
obvious way, so we can view it (\scite{\FD}{VIII, 3.8}) as a central
\cstar-subalgebra of the multiplier algebra of $C^*(\D)$.

Let $J_{\hbar}=\{f\in C_0(I): f(\hbar)=0\}$. It only remains to show
that $C^*(\CE^{\hbar})=C^*(\D)J_{\hbar}$. For then, by
\scite{\RieffelCF}{1.2}, we will have that $\hbar\mapsto
\|\tilde{\Pi}^{\hbar}(\phi)\|$ is upper semicontinuous for all $\phi\in
C^*(\D)$. This implies that $\hbar\mapsto \|\psi\|_{\hbar}$ is upper
semicontinuous for any $\psi\in \ccb$.  Now, it is apparent that $\phi
j\in L^1(\CE)$ for $j\in J_{\hbar}$, and $\phi\in L^1(\D)$, which shows
that $C^*(\D)J_{\hbar}\subset C^*(\CE^{\hbar})$.  On the other hand, if
$\{e_{\lambda}\}$ is a bounded approximate identity for $J_{\hbar}$,
then $\lim_{\lambda} \phi e_{\lambda}=\phi$ for all $\phi\in
C^*(\CE^{\hbar})$: It suffices to show it for compactly supported maps
$\phi$, since $\{e_{\lambda}\}$ is assumed to be bounded.  Notice that
the statement holds for $\phi=f\delta_e$, with $f\in E_e^{\hbar}$,
because $E_e^{\hbar}\cong B_e\otimes J_{\hbar}$. Now, if
$\phi=f_x\delta_x$ for some $f_x\in E^{\hbar}_x$, we have
  $$
  \|\phi e_{\lambda}-\phi\|^2=\|(\phi e_{\lambda}-\phi)^*(\phi
e_{\lambda} -\phi)\|\leq (\|e_{\lambda}\|+1)\|\phi^*\phi e_{\lambda}-
\phi^*\phi\|,
  $$
  which goes to zero because $\phi^*\phi \in E^{\hbar}_e$. This shows
that $C^*(\D)J_{\hbar}\supset C^*(\CE^{\hbar})$.
  \proofend

  \state Proposition If $G$ is also amenable, then the map $\hbar\mapsto
\|\phi\|_{\hbar}$ is lower semicontinuous on $I$ for all $\phi\in \ccb$.

  \proof Since $G$ is amenable, the left regular representation
$\Lambda^{\hbar}$ of $C^*(\B^{\hbar})$ is faithful (\scite{\Amena}{2.3
and 4.7}), so it suffices to show that $\hbar\mapsto
\|\Lambda^{\hbar}_{\phi}\|$ is lower semicontinuous for $\phi\in \ccb$.

As in \cite{\Amena}, for $\hbar\in I$ we denote by $L^2(\B^{\hbar})$ the
completion of $C_c(\B^{\hbar})$ with its obvious right pre-Hilbert
module structure over $B_e^{\hbar}$, which yields the norm
  $$
  \|\xi\|^2 =
  \|\sum_{x\in G} \xi(x)\starh \prodh \xi(x)\|_{B_e^{\hbar}} =
  \|\sum_x\theta^{\hbar}_{x^{-1}}[\xi(x)^*\xi(x)]\|_{B^0_e},
  $$
  for any $\xi\in \ccb$, the undecorated involution and multiplication
denoting those in $\B^0$.

The left regular representation $\Lambda^{\hbar}$ of $\phi\in \ccb$ is
the adjointable operator given by:
  $$
  (\Lambda^{\hbar}_{\phi}\xi)(y) =
  \sum_{x\in G} \phi(x)\prodh \xi(x^{-1}y) =
  \sum_x\phi(x)\theta^{\hbar}_x[\xi(x^{-1}y)],
  $$
  for $\xi\in \ccb\subset L^2(\B^{\hbar})$. So we have
  $$
  \|\Lambda^{\hbar}_{\phi}\xi\|^2_{\hbar} =
  \|\sum_{x,y}\theta^{\hbar}_{y^{-1}}[(\phi(x)\theta^{\hbar}_x
    (\xi(x^{-1}y)))^*
    (\phi(x)\theta^{\hbar}_x(\xi(x^{-1}y)))]\|_{B_e}.
  $$

Notice that the sum above is finite, since both $\phi$ and $\xi$ are
compactly supported. Besides, each term of the sum is continuous on
$\hbar$, so $\hbar\mapsto \|\Lambda^{\hbar}_{\phi}\xi\|_{\hbar}$ is
continuous. Now fix $\phi\in \ccb$, $\epsilon >0$, and $\hbar_0\in
I$. Then $\xi_0\in C_c(\B^{\hbar_0})$ can be found so that $\|x_0\|=1$
and
$\|\Lambda^{\hbar_0}_{\phi}\xi_0\|>\|\Lambda^{\hbar_0}_{\phi}\|-\epsilon$. For
one can find $\xi\in L^2(\B^{\hbar_0})$ satisfying that inequality for
$\frac{\epsilon}{2}$, with $\|\xi\|=1$. Then, given $\{\xi_n\}\subset
\ccb$ such that $ \lim \xi_n= \xi$, the sequence
$\{\frac{1}{\|\xi_n\|}\xi_n\}$ also converges to $\xi$. So one can take
$\xi_0\in \ccb$, such that $\|\xi_0\|=1$ and $\|\xi-\xi_0\|<
\frac{\epsilon}{2}\|\Lambda^{\hbar_0}_{\phi}\|$. Then
  $$ \|\Lambda^{\hbar_0}_{\phi}\|-\frac{\epsilon}{2}<
\|\Lambda^{\hbar_0}_{\phi}\xi_0\| + \frac{\epsilon}{2},
  $$
  as required. It now follows that, for $\hbar$ close enough to
$\hbar_0$,
  $$
  \frac{\|\Lambda^{\hbar}_{\phi}\xi_0\|_{\hbar}}{\|\xi_0\|_{\hbar}}>
\|\Lambda^{\hbar_0}_{\phi}\|-\epsilon \hbox{, so }
\|\Lambda^{\hbar}_{\phi}\|> \|\Lambda^{\hbar_0}_{\phi}\|-\epsilon.
  $$
  \proofend

We summarize the previous results in the following theorem.

  \state \label\ContThm Theorem Let $\B$ be a Fell bundle over a
discrete amenable group $G$, and let $B=C^*(\B)$. If
$\{\theta^{\hbar}:\hbar\in I\}$ and $\B^{\hbar}$ are as in
\lcite{\defi}, then $\{C^*(\B^{\hbar}), \Lambda\}$ is a continuos field
of \cstar-algebras, such that $C^*(\B^0)=B$, where $\Lambda$ is the
family of cross-sections obtained, as in \scite{\Dixmier}{10.2.3}, out
of $C_c(\B^{\hbar})$.

  \endgroup
  \section {Discrete abelian groups}
  We would now like to describe a method for producing examples of the
above situation.  To reduce the technical difficulties to a minimum we
will consider here exclusively the case of discrete abelian groups.
Several interesting examples, however, will fit this context.

  Fix, throughout this section, a discrete abelian group $G$ and let
$\Gh$ be its Pontryagin dual, so that $\Gh$ is a compact abelian group.
We shall denote the duality between $G$ and $\Gh$ by
  $$
  (x,t) \in \Gh \x G \vaipara \<x,t\> \in S^1.
  $$

  Let $B$ be a \cstar-algebra carrying a continuous action $\g$ of
$\Gh$.  For each $t$ in $G$, the $t$-spectral subspace of $B$ is defined
by
  $$
  B_t = \{ b\in B : \g_x(b) = \<x,t\> b, \hbox{ for all } x \in \Gh\}.
  $$

  It is an easy exercise to show that each $B_t$ is a closed linear
subspace of $B$, that $B_r B_s \subseteq B_{rs}$, and that $B_t^* =
B_{t\inv}$.  By imitating \scite{\Circle}{2.5} one can show that $B$
coincides with the closure of
  $\bigoplus_{t \in G} B_t$
  (we use the symbol $\bigoplus$ to denote the algebraic direct sum,
that is, the set of \stress{finite} sums)
  and that the formula
  $$
  P_t (b) = \int_{\Gh} \<x,t\>\inv \g_x (b)\, dx \for b\in B, t\in G,
  $$
  defines a contractive projection $P_t$, from $B$ onto $B_t$, where the
integral is taken with respect to normalized Haar measure on $\Gh$.  If
$e$ denotes the unit of $G$, then $P_e$ is in fact a positive
conditional expectation onto $B_e$.

The collection $\B = \{B_t\}_{t\in G}$ therefore constitutes a Fell
bundle over $G$ and also makes $B$ into a topologically $G$-graded
algebra, as defined in \scite{\Amena}{3.4}.  If we take into account the
fact that abelian groups are amenable, and use
  \scite{\Amena}{4.7}
  in combination with
  \scite{\Amena}{4.2},
  we then conclude that $B$ is isomorphic to the (full) cross-sectional
\cstar-algebra of $\B$
  \scite{\FD}{VIII.17.2}
  as well as to its reduced cross-sectional \cstar-algebra
\scite{\Amena}{2.3}.

  Now suppose that, in addition to the action $\g$ above, we are given
an action $\t$ of $G$ on $B$ which commutes with $\g$, in the sense that
each $\g_x$ commutes with each $\t_t$.  It then follows that $\t_s(B_t)
\subseteq B_t$ for each $t,s$ in $G$, so that we can think of $\t$ as an
action of $G$ on the Fell bundle $\B$.

  This can in turn be fed to the construction described in the previous
section, providing the $\t$-deformed bundle $\DBs$.  Taking a further
step, one can form the cross-sectional algebra of this bundle.

  \definition Given commuting actions $\g$ and $\t$, respectively of
$\Gh$ and $G$, on the \cstar-algebra $B$, the cross-sectional algebra of
$\DBs$ will be called the $(\g,\t)$-deformation of $B$ and will be
denoted $\DAs$.

  It should be noted that, if $\t$ is the trivial action, then $\DBs$ is
nothing but $\B$ itself and hence, by the comment made earlier, its
cross-sectional \cstar-algebra coincides with $B$, and so $\DAs = B$.
Likewise, if $\g$ is trivial then $B_t = \{0\}$, for all $t$, except for
$B_e$ which is the whole of $B$ and, once more, one has that $\DAs = B$.
However, if neither group acts trivially, then the algebraic structure
of $B$ may suffer a significant transformation as it will become
apparent after we discuss a few examples.

  We like to think of this as if the pair $(\g,\t)$ ``causes'' a
deformation on $B$.  This motivates our next:

  \definition A \stress{\dd{}} for a \cstar-algebra $B$ consists of a
triple
  $(G,\g,\t)$,
  where $G$ is a discrete abelian group, and $\g$ and $\t$ are commuting
actions, respectively of $\Gh$ and $G$, on $B$.  The action $\g$ will be
called the
  \stress{gauge}
  action while $\t$ will be referred to as the
  \stress{deforming} action.

  Unless otherwise noted, whenever we speak of the Fell bundle
  $\B = \{B_t\}_{t\in G}$, in the presence of a \dd $(G,\g,\t)$ for a
\cstar-algebra $B$, we will be referring to the spectral decomposition
of the gauge action, as above.

  \sysstate
  {Remark}
  {\rm}
  {\label{\Remark} Observe that $\DAs$, being the cross-sectional
\cstar-algebra of $\DBs$, contains the algebraic direct sum
  $\bigoplus_{t \in G} B_t$
  as a dense *-sub-algebra.
  It should be noted that the set
  $\bigoplus_{t \in G} B_t$
  itself, as well as its linear structure, depends exclusively on the
gauge action.  However, its involution and multiplication operations are
strongly dependent of the deforming action.
  Also,
  since the process of deformation does not affect the norm structure of
the Fell bundle, and since
  the fibers of that bundle embed isometricaly into its cross-sectional
algebra, we see that the norm of an element belonging to a single fiber
remains unaffected by the deformation.  However, there is not much we
can say about the norm of other elements in $\bigoplus_{t \in G} B_t$.
  Summarizing, in case we are given several \dds sharing the same gauge
action, it will be convenient to think of the deformed algebras as
completions of
  $\bigoplus_{t \in G} B_t$
  under different norms and with different algebraic operations.
  }

  \state Proposition \label{\TildeAction} Let $(G,\g,\t)$ be a \dd for a
\cstar-algebra $B$.  Suppose $B$ carries a third continuous action
$\alpha$, this time of a locally compact group $H$, which commutes both
with $\g$ and $\t$.  Then there exists a continuous action
$\tilde\alpha$ of $H$ on $\DAs$ which coincides with $\alpha$ on
  $\bigoplus_{t \in G} B_t$.

  \proof Since $\alpha$ commutes with the gauge action, each spectral
subspace $B_t$ is invariant by $\alpha_h$, for each $h\in H$.  So
$\alpha_h$ can be thought of as an automorphism of the Fell bundle $\B$.
We claim it is also automorphic for the deformed structure.  In fact, if
  $b_t\in B_t$ and $b_s\in B_s$ then
  $$
  \alpha_h(b_t \defprod b_s) =
  \alpha_h(b_t \t_t(b_s)) =
  \alpha_h(b_t) \t_t(\alpha_h(b_s)) =
  \alpha_h(b_t) \defprod \alpha_h(b_s),
  $$
  and
  $$
  \alpha_h(b_t\defstar) =
  \alpha_h(\t_t\inv(b_t^*)) =
  \t_t\inv(\alpha_h(b_t)^*) =
  \alpha_h(b_t)\defstar.
  $$

  Thus $\alpha_h$ extends to an automorphism of $\DAs$.  The remaining
verifications are left to the reader.
  \proofend

  Among the possible choices for the action $\alpha$ above one could
take the gauge action itself, so one can speak of the ``deformed gauge
action'', that is $\tilde\g$.

  \state Proposition \label{\spectilde} For each $t$ in $G$, the
$t$-spectral subspace for the deformed gauge action on $\DAs$ is
precisely $B_t$.

  \proof Let us temporarily denote the $t$-spectral subspace for
$\tilde\g$ by $\tilde B_t$.  Since $\tilde\g$ coincides with $\g$ on
  $\bigoplus_{t \in G} B_t$,
  it is clear that
  $
  \tilde\g_x(b_t) = \<x,t\> b_t
  $
  for each $b_t$ in $B_t$.  So $B_t \subseteq \tilde B_t$.  Conversely,
let $a\in \tilde B_t$.  Then, for each $\varepsilon > 0$, take a finite
sum
  $\sum_{r\in G} b_r$
  with $b_r \in B_r$, and such that
  $\| a - \sum_{r\in G} b_r \| < \varepsilon$.
  Considering the spectral projections
  $$
  \tilde P_t (b) =
  \int_{\Gh} \<x,t\>\inv \tilde\g_x (b)\, dx \for b\in \DAs, t\in G,
  $$
  we have $a=\tilde P_t(a)$ while
  $
  \tilde P_t(\sum_{r\in G} b_r) = b_t
  $.
  So
  $
  \| a - b_t \| =
  \| \tilde P_t(a - \sum_{r\in G} b_r) \|
  < \varepsilon
  $.
  This says that $a$ is in the closure of $B_t$ within $\DAs$.  But
since the norm on $B_t$ is not affected by the deformation, it remains a
Banach space after the deformation is performed, and hence it is closed
in $\DAs$.  Therefore $a\in B_t$.
  \proofend

  \state Theorem \label{\FixedPoint} Let $(G,\g,\t)$ be a \dd for a
\cstar-algebra $B$, and let $\alpha$ be an action of a group $H$ on $B$
which commutes both with $\g$ and $\t$.  Let $B^0$ be the fixed point
sub-algebra of $B$ for $\alpha$, and let $\g^0$ and $\t^0$ be the
restrictions of $\g$ and $\t$ to $B^0$, respectively.  Then the deformed
algebra $\DA{(B^0)}{\g^0}{\t^0}$ is isomorphic, in a natural way, to the
fixed point sub-algebra of $\DAs$ for $\tilde\alpha$.

  \proof Observe, initially, that since $\tilde\alpha$ and $\tilde\g$
coincide with $\g$ and $\alpha$, respectively, on
  $\bigoplus_{t \in G} B_t$,
  then they must commute.

  This implies that the fixed point sub-algebra for $\tilde\alpha$,
which we denote by $A$, is invariant under $\tilde\g$.  Because of
\lcite{\spectilde}, it is a simple matter to verify that the spectral
decomposition of the restriction of $\tilde\g$ to $A$ is given by
  $
  \bigoplus_{t \in G} B_t \cap A
  $.
  Observe that, because $\alpha$ and $\tilde\alpha$ agree on each $B_t$,
  $$
  B_t \cap A =
  \{b\in B_t : \alpha_h(b) = b \hbox{ for all } h\in H\}
  \$=
  \{b\in B : \alpha_h(b) = b
              \hbox{ for all } h\in H
              \hbox{ and }
              \g_x(b) = \<x,t\> b
              \hbox{ for all } x \in \Gh\}
  \$=
  \{b\in B^0 : \g_x(b) = \<x,t\> b, \hbox{ for all } x \in \Gh\} =
  B_t^0,
  $$
  where we have denoted by $B_t^0$ the $t$-spectral subspace of $B^0$
under $\g^0$.  It is now easy to see that the Fell bundle structure
arising from the grading $\{B_t \cap A\}_{t\in G}$ of $A$, and that of
the grading of the deformed algebra
  $\DA{(B^0)}{\g^0}{\t^0}$ are isomorphic.  The result then follows from
  \scite{\Amena}{4.2}.
  \proofend

  \section {The derivative of the deformed product}
  \global\edef\DerivSection{\number\secno}
  Let $B$ be a \cstar-algebra carrying a strongly continuous action
$\action$ of $\R^{2d}$.

  For each $j=1,\ldots,2d$, define the differential operator
$\partial_{u_j}$ on $B$ by
  $$
  \pdop{u_j}{f}
  =
  \deriv{\action_{(0,\ldots,\lambda,\ldots,0)}(f)}{\lambda}{0}
  \for f\in B,
  $$
  where the $\lambda$ in $(0,\ldots,\lambda,\ldots,0)$ appears in the
$j^{th}$ position.  Of course $\pdop{u_j}{f}$ is only defined when $f$
is sufficiently smooth.  In particular this is the case for the
\stress{$\action$-smooth} elements, that is, those elements $f\in B$
such that
  $$
  u\in \R^{2d} \vaipara \action_u(f) \in B
  $$
  is an infinitely differentiable Banach space valued function.  It is
well known that these elements form a dense subset of $B$ (see, e.g,
\scite{\Brat}{2.2.1}).

  In what follows we shall adopt the coordinate system
$(x_1,\ldots,x_d,y_1,\ldots,y_d)$ on $\R^{2d}$ and hence we shall speak
of the differential operators $\partial_{x_j}$ and $\partial_{y_j}$, for
$j=1,\ldots,d$.

  In \cite{\RieffelDQi} (see also \cite{\RieffelDQii}) Rieffel showed
how to construct a strict deformation quantization of $B$ ``in the
direction'' of the Poisson bracket $\{\cdot,\cdot\}$ defined by
  $$
  \{f,g\} = \sum_{j=1}^d \pdop{x_j}{f} \pdop{y_j}{g} - \pdop{y_j}{f}
\pdop{x_j}{g},
  $$
  in the important special case when $B$ is the algebra of continuous
functions on a smooth manifold.  Rieffel deals, in fact, with a more
general situation, where the Poisson bracket involves the choice of a
skew-symmetric matrix $J$.

  Without attempting to develop a general theory, we would now like to
describe a connection between Rieffel's theory and ours.  Our goal will
be to compute the derivative of the deformed product on $B$, arising
{}from a certain \dd associated to $\action$.
  The technical complications will be kept to a minimum by assuming that
$\action$ is periodic in the first $d$ variables.

  Let $\g$ be the action of $\R^d$ given by the restriction of $\action$
to its first $d$ variables, that is
  $$
  \g_{(x_1,\ldots,x_d)} = \action_{(x_1,\ldots,x_d,0,\ldots,0)}
  \for (x_1,\ldots,x_d) \in \R^d.
  $$
  Our periodicity assumption, to make it precise, is that $\g$ is
trivial on $\Z^d$ and hence it defines, by passage to the quotient
$\R^d/\Z^d$, an action of the $d$ dimensional torus $\T^d$ on $B$, which
we will still denote by $\g$.

  On the other hand, consider the action $\t$ of $\R^d$ on $B$ defined
by
  $$
  \t_{(y_1,\ldots,y_d)} =
  \action_{(0,\ldots,0,y_1,\ldots,y_d)}
  \for (y_1,\ldots,y_d) \in \R^d.
  $$
  If $\hbar$ is a real number, we will let the action $\t^\hbar$ of
  $\Z^d$ on $B$ be defined by
  $$
  \t^\hbar_{(n_1,\ldots,n_d)} = \t_{(\hbar n_1,\ldots,\hbar n_d)}
  \for (n_1,\ldots,n_d) \in \Z^d.
  $$

  Since both $\g$ and $\t^\hbar$ come from the action of the commutative
group $\R^{2d}$, it is clear that they commute with each other.
Moreover, since the Pontryagin dual of the group $\Z^d$ is precisely
$\T^d$, the triple
  $(\Z^d,\g,\t^\hbar)$
  is seen to be a \dd for $B$.

  Let $\B = \{B_t\}_{t\in G}$ be the Fell bundle arising from the
spectral decomposition of $\gamma$.  We may then speak of the deformed
bundle $\DB{\B}{\t^\hbar}$ whose operations will be denoted by
$\defprod_\hbar$ and $\null^{\defstarsymbol_\hbar}$.  We also have the
deformed algebra $\DA{B}{\g}{\t^\hbar}$, which we will simply denote by
$B^{(\hbar)}$.

  \state Proposition If $f$ is $\action$-smooth then $P_t(f)$ is also
$\action$-smooth for all $t$ in $\Z^d$.
  In addition, for $j=1,\ldots,2d$, we have
  $
  \pdop{u_j}{P_t(f)} =
  P_t(\pdop{u_j}{f}),
  $
  and therefore each $B_t$ is invariant under $\partial_{u_j}$.

  \proof For $u\in \R^{2d}$ we have
  $$
  \action_u(P_t(f))
  =
  \action_u
    \( \int_{\T^d} \<x,t\>\inv \g_x (f)\,dx \)
  =
  \int_{\T^d} \<x,t\>\inv \g_x (\action_u(f))\,dx,
  $$
  which is therefore smooth as a function of $u$.  This shows that
$P_t(f)$ is $\action$-smooth.  We have
  $$
  \pdop{u_j}{P_t(f)}
  =
  \deriv{\action_{(0,\ldots,\lambda,\ldots,0)}P_t(f)}{\lambda}{0}
  \$=
  \int_{\T^d}\deriv{
    \action_{(0,\ldots,\lambda,\ldots,0)}
      \( \<x,t\>\inv \g_x (f) \)}
    {\lambda}{0}\,dx
  \$=
  \int_{\T^d} \<x,t\>\inv \g_x (\pdop{u_j}{f})\,dx
  =
  P_t( \pdop{u_j}{f} ).
  \proofend
  $$

  \state Lemma \label{\Single} Let
  $t=(t_1,\ldots,t_d)$
  and
  $s=(s_1,\ldots,s_d)$
  be in $\Z^d$ and take $f\in B_t$ and $g \in B_s$.  Suppose that $g$ is
smooth for $\t$.  Then, for all real numbers $\hbar$
  $$
  \left\Vert {f \defprod_\hbar g - fg \over \hbar} -
  {1 \over 2\pi i}
  \sum_{j=1}^d \pdop{x_j}{f} \pdop{y_j}{g} \right\Vert
  \leq
  |\hbar|\,\|f\|
  \left\Vert
  \sum_{j,k=1}^d t_j t_k \pdop{y_j}{\pdop{y_k}{g}}
  \right\Vert.
  $$

  \proof Initially we would like to stress that the term whose norm is
referred to, in the left hand side above, lies in $B_{t+s}$.  This
Banach space embeds isometricaly into each $B^{(\hbar)}$, and hence its
norm is unambiguously defined.  We have
  $$
  f \defprod_\hbar g - fg
  = f\t^\hbar_t(g) - fg.
  $$
  Now, consider the $C^\infty$ map
  $
  F : \R \arw B
  $
  given by
  $$
  F(\hbar) \=
  \t^\hbar_t(g) =
  \action_{(0,\ldots,0,\hbar t_1,\ldots,\hbar t_d)}(g).
  $$
  Its first two derivatives are given by
  $$
  F'(\hbar) = \action_{(0,\ldots,0,\hbar t_1,\ldots,\hbar t_d)}
    \( \sum_{j=1}^d t_j \pdop{y_j}{g} \),
  $$
  and
  $$
  F''(\hbar) = \action_{(0,\ldots,0,\hbar t_1,\ldots,\hbar t_d)}
    \( \sum_{j,k=1}^d t_j t_k\pdop{y_j}{\pdop{y_k}{g}} \),
  $$
  for all $\hbar$ in $\R$.
  The first order Taylor expansion for $F$ reads
  $$
  F(\hbar) = F(0) + \hbar F'(0) +
    \int_0^\hbar (\hbar-\lambda)F''(\lambda)\,d\lambda,
  $$
  from where we conclude that
  $$
  \left\Vert
  {F(\hbar) - F(0) \over \hbar} - F'(0)
  \right\Vert
  \leq
  |\hbar| \sup_{\lambda \in I} \| F''(\lambda) \|,
  $$
  where the interval $I$ is either $[0,\hbar]$ or $[\hbar,0]$, depending
on the sign of $\hbar$.  Translating this back in terms of $g$, we
conclude that
  $$
  \left\Vert
  {\t^\hbar_t(g) - g \over \hbar}
  -
  \sum_{j=1}^d t_j \pdop{y_j}{g}
  \right\Vert
  \leq
  |\hbar|
  \left\Vert
  \sum_{j,k=1}^d t_j t_k\pdop{y_j}{\pdop{y_k}{g}}
  \right\Vert.
  $$
  Using the first equation obtained in the course of the present proof
gives
  $$
  \left\Vert
  {f \defprod_\hbar g - fg \over \hbar}
  -
  \sum_{j=1}^d t_j f \pdop{y_j}{g}
  \right\Vert
  \leq
  |\hbar|\,\|f\|
  \left\Vert
  \sum_{j,k=1}^d t_j t_k\pdop{y_j}{\pdop{y_k}{g}}
  \right\Vert.
  $$
  On the other hand, recall that $f$ is in the $t$-spectral subspace of
the gauge action.  This means that, for
  $x = (x_1,\ldots,x_d) \in \R^d$,
  we have that
  $\g_x(f) = \<x,t\>f$,
  or
  $$
  \g_x(f) = \expi{x_1 t_1} \ldots \expi{x_d t_d}f.
  $$
  If follows that
  $\pdop{x_j}{f} = 2\pi i t_j f$, and hence that
  $t_j f = (2\pi i)\inv\pdop{x_j}{f}$, which, when plugged into the last
inequality above, leads to the conclusion.
  \proofend

  The purpose of this Lemma, as the reader may have anticipated, is to
allow us to compute the derivative of
  $f \defprod_\hbar g$,
  with respect to $\hbar$, which is one of the most important
ingredients in Rieffel's theory of deformation quantization
\cite{\RieffelDQi, \RieffelDQii}.  However the expression
  $f \defprod_\hbar g$, strictly speaking, applies only for $f$ and $g$
belonging, each, to a single spectral subspace of the gauge action.  The
question we want to address is this:

  \state Question \label{\Question} What is the biggest subset of $B$
which can be mapped, in a natural way, into each deformed algebra
$B^{(\hbar)}$?

  The remark made in \lcite{\Remark} is relevant here, providing
  $\bigoplus_{t \in \Z^d} B_t$
  as a partial answer.  Pushing this further, recall that the
cross-sectional \cstar-algebra of our Fell bundle $\B$ is defined to be
the enveloping \cstar-algebra of the $L_1$ cross-sectional algebra
$L_1(\B)$.  Therefore each $B^{(\hbar)}$ contains a copy of $L_1(\B)$
which, again by \lcite{\Remark}, does not depend on $\hbar$, as far as
its normed linear space structure is concerned.
  A better answer to our question is thus $L_1(\B)$.

  We do not claim, however, that this is the best possible answer.  In
fact, the word \stress{natural} in \lcite{\Question} lacks a precise
meaning, as it stands.  The correct way to rephrase \lcite{\Question}
could possibly be:

  \state Question For each $\hbar$, let
  $\iota^\hbar : L_1(\B) \arw B^{(\hbar)}$
  be the natural inclusion, viewed as a densely defined linear map on
$B$.  Is $\iota^\hbar$ closable?  That is, is the closure of its graph,
the graph of a well defined linear map?  If so, how to characterize the
domain $D^\hbar$ of this map?  Is there any relationship between the
$D^\hbar$ for different $\hbar$?  What is the intersection of the
$D^\hbar$ as $\hbar$ ranges in $\R$?

  In defense of the $L_1$ cross-sectional algebra we must say that it
includes the smooth elements for the gauge action: it is a well known
fact that, for such an element $f$, one has that
  $f = \sum_{t\in \Z^d} P_t(f)$,
  where the series is absolutely convergent, since it satisfies
Schwartz's condition, namely that
  $$
  \sup_{t\in \Z^d} \|h(t) P_t(f) \| < \infty,
  $$
  for every complex polynomial $h$ in the $d$ variables
$t_1,\ldots,t_d$.

  \state Theorem \label{\Derivative} Let $f,g\in B$ be $\action$-smooth
elements.  Then
  $$
  \lim_{\hbar\rightarrow 0}
  \left\Vert {f \defprod_\hbar g - fg \over \hbar} -
  {1 \over 2\pi i}
  \sum_{j=1}^d \pdop{x_j}{f} \pdop{y_j}{g} \right\Vert_\hbar = 0,
  $$
  where $\|\cdot\|_\hbar$ refers to the norm of the deformed algebra
$B^{(\hbar)}$.

  \proof Initially we should observe that the terms appearing between
the double bars, above, all have natural interpretations as elements of
$B^{(\hbar)}$.  This is because the smooth elements
  $f$,
  $g$,
  $fg$, and
  $\pdop{x_j}{f} \pdop{y_j}{g}$,
  may be seen as elements of
  $L_1(\B)$, which, in turn, may be interpreted as a subset of
$B^{(\hbar)}$, according to the comment above.

  Write
  $f = \sum_{t\in \Z^d} P_t(f)$ and
  $g = \sum_{t\in \Z^d} P_t(g)$.  For each $j=1,\ldots,2d$ we have that
$\pdop{u_j}{f}$ is also smooth, hence it ``Fourier series'' converges:
  $$
  \pdop{u_j}{f}
  =
  \sum_{t\in \Z^d} P_t(\pdop{u_j}{f})
  =
  \sum_{t\in \Z^d} \pdop{u_j}{P_t(f)},
  $$
  and similarly for $g$.  So,
  $$
  \sum_{j=1}^d \pdop{x_j}{f} \pdop{y_j}{g}
  =
  \sum_{t,s\in \Z^d}
    \sum_{j=1}^d \pdop{x_j}{P_t(f)} \pdop{y_j}{P_s(g)}.
  $$
  Also
  $$
  {f \defprod_\hbar g - fg \over \hbar}
  =
  \sum_{t,s\in \Z^d}
    {P_t(f) \defprod_\hbar P_s(g) - P_t(f)P_s(g) \over \hbar}.
  $$
  Using \lcite{\Single}, it follows that
  $$
  \left\Vert {f \defprod_\hbar g - fg \over \hbar} -
  {1 \over 2\pi i}
  \sum_{j=1}^d \pdop{x_j}{f} \pdop{y_j}{g} \right\Vert_\hbar
  \$\leq
  \sum_{t,s\in \Z^d}
  \left\Vert {P_t(f) \defprod_\hbar P_s(g) - P_t(f)P_s(g) \over \hbar} -
  {1 \over 2\pi i}
  \sum_{j=1}^d \pdop{x_j}{P_t(f)} \pdop{y_j}{P_s(g)} \right\Vert_\hbar
  \$\leq
  |\hbar|
  \sum_{t,s\in \Z^d}
  \|P_t(f)\|
  \left\Vert
  \sum_{j,k=1}^d t_j t_k \pdop{y_j}{\pdop{y_k}{P_s(g)}}
  \right\Vert
  \$\leq
  |\hbar|
  \sum_{j,k=1}^d
    \( \sum_{t\in \Z^d} |t_j t_k|\,\|P_t(f)\| \)
    \( \sum_{s\in \Z^d} \|
      P_s(\pdop{y_j}{\pdop{y_k}{g}})
    \| \).
  $$
  By our hypothesis, these infinite series converge, and hence the whole
thing tends to zero as $\hbar \rightarrow 0$.
  \proofend

  \sysstate{Remark}{\rm}
  {If one is interested in determining the exact class of
differentiability needed for the above result to hold, a quick look at
the last displayed expression, in the proof above, gives the answer.
That is, $f$ should be supposed to be of class
  $C^{2d+2}$ for $\g$, and the second order differential of $g$ with
respect to $\t$ should be of class
  $C^{2d}$ for $\g$.  These conditions imply the convergence of these
infinite series, and hence the conclusion.}

  Our next result shows that the infinitesimal commutator, for the
deformed product, is given by the Poisson bracket described at the
beginning of this section.  Its proof is an immediate consequence of
\lcite{\Derivative}.

  \state Corollary \label{\PoissonCorollary} Let $f,g\in B$ be smooth
elements for $\action$.  Then
  $$
  \lim_{\hbar\rightarrow 0}
  \left\Vert {f \defprod_\hbar g - g \defprod_\hbar f - [f,g]
  \over \hbar} -
  {1 \over 2\pi i}
  \{f,g\}
  \right\Vert_\hbar = 0,
  $$
  where
  $[\cdot,\cdot]$ is the commutator for the original multiplication on
$B$,
  and
  $\{\cdot,\cdot\}$ is the Poisson bracket defined near the begining of
this section.

  Since the family
  $\{\t^\hbar\}_{\hbar\in \R}$
  is obviously continuous in the sense of section
  \lcite{\ContSection},
  we get, by
  \lcite{\ContThm},
  a continuous field of \cstar-algebras
  $\{B^{(\hbar)}\}_{\hbar\in \R}$, and hence a strict deformation
quantization in the sense of Rieffel
  \scite{\RieffelHeis}{Definition 1.1},
  with the modification, required in the noncommutative situation,
  corresponding to the introduction of the term $[f,g]$ in the statement
of \lcite{\PoissonCorollary}.  With this remark we have, for future
reference, the following:

  \state Corollary The family
  $\{B^{(\hbar)}\}_{\hbar\in \R}$
  gives a strict deformation quantization for $B$, in the direction of
the Poisson bracket defined above.

  \section {Example: Non commutative 3-spheres}
  In \cite{\Mt} Matsumoto defined a family of \cstar-algebras, denoted
  $S^3_\tt$,
  depending on a real parameter $\tt$.  This family is to be thought of
as a deformation of the commutative \cstar-algebra $C(S^3)$ of all
continuous complex valued functions on the 3-sphere $S^3$, because, when
$\tt=0$, $S^3_\tt$ is isomorphic to $C(S^3)$.

The purpose of the present section is to show that $S^3_\tt$ can be
constructed from a certain \dd for the algebra $C(S^3)$.  Recall from
\cite{\Mt} that $S^3_\tt$ may be defined as the universal \cstar-algebra
given by generators and relations as follows: for generators take
symbols $S$ and $T$ and for relations consider
  \smallskip\itemitem{M-1)} $S^*S = S S^*$, $T^*T = T T^*$,
  \smallskip\itemitem{M-2)} $\|S\| \leq 1$, $\|T\| \leq 1$,
  \smallskip\itemitem{M-3)} $(1-T^*T) (1-S^*S) = 0$, and
  \smallskip\itemitem{M-4)} $TS = \expi{\tt} ST$.

  \medskip An alternative description of $S^3_\tt$ is given by
\scite{\Mt}{8.1}.  It says that $S^3_\tt$ is also the universal
\cstar-algebra on the generators $B$ and $C$ satisfying
  \smallskip\itemitem{M-1')} $B^*B = B B^*$, $C^*C = C C^*$,
  \smallskip\itemitem{M-2')} $B^*B + C^*C = 1$, and
  \smallskip\itemitem{M-3')} $CB = \expi{\tt} BC$.

  \medskip The relationship between these presentations is given by the
formulas
  $$
  B = S(S^*S + T^*T)^{-{1 \over 2}}, \quad
  C = T(S^*S + T^*T)^{-{1 \over 2}}.
  $$

  Define an action $\g$ of $S^1$ on $S^3$ by
  $$
  \g_\lambda(z,w) = (\lambda z, \lambda w),
  $$
  where $z,w, \lambda\in \C$ satisfy $|z|^2 + |w|^2 = 1$ and
$|\lambda|=1$.

  Also, fixing a real number $\tt$, define an action $\t$ of $\Z$ on
$S^3$ by
  $$
  \t_n(z,w) = (\expi{n\tt}z,w) \for (z,w)\in S^3, \quad n\in \Z.
  $$
  These give actions of $S^1$ and $\Z$ on $C(S^3)$ by letting
  $$
  \g_\lambda(f)\calcat{(z,w)} = f(\lambda z, \lambda w), \and
  \t_n(f)\calcat{(z,w)} = f(\expi{n\tt}z,w)
  $$
  for $f\in C(S^3)$, $(z,w)\in S^3$, $\lambda \in S^1$ and $n\in\Z$.
Noting that $\g$ and $\t$ commute with each other, we see that we are
facing a \dd $(\Z,\g,\t)$ for the algebra $C(S^3)$.

  \state Theorem \label{\NCSphere} The deformed algebra
  $\DA{C(S^3)}{\g}{\t}$
  is isomorphic to Matsumoto's algebra $S^3_\tt$.

  \proof Let $Z,W\in C(S^3)$ be the functions defined by
  $$
  Z(z,w) = z, \and
  W(z,w) = w,
  $$
  for $(z,w)\in S^3$.
  Since $\g_\lambda(Z) = \lambda Z$ and $\g_\lambda(W) = \lambda W$, we
have that both $Z$ and $W$ belong to the first spectral subspace for
$\g$.  Then, regarding the deformed product, we have
  $$
  Z\defprod W = Z \t_1(W) = ZW
  \and
  W\defprod Z = W \t_1(Z) = \expi{\tt} WZ,
  $$
  so that
  $$
  W\defprod Z = \expi{\tt}Z\defprod W.
  $$
  This says that $Z$ and $W$ satisfy (M-3').  It is easy to check that
they also satisfy (M-1') and (M-2') with respect to the deformed product
an involution.  So, by the universal property, there exists a
\cstar-algebra homomorphism
  $$
  \psi : S^3_\tt \arw
  \DA{C(S^3)}{\g}{\t},
  $$
  such that
  $\psi(B) = Z$ and $\psi(C) = W$, which we claim to be an isomorphism.

  To show that $\psi$ is surjective observe that, since $Z$ and $W$
belong to the image of $\psi$, we just have to show that $Z$ and $W$
generate
  $\DA{C(S^3)}{\g}{\t}$.
  In the special case of $\tt = 0$ this of course follows from the
Stone--Weierstrass Theorem.  Looking closer we can actually show that
the $n$-spectral subspace for the action $\g$ on $C(S^3)$ is linearly
spanned by the set
  $$
  \{Z^i {Z^*}^j W^k {W^*}^l : i,j,k,l \in {\bf N}, i-j+k-l=n \}.
  $$

  In fact, any $f\in C(S^3)$ may be arbitrarily approximated by a linear
combination of terms of the form $Z^i {Z^*}^j W^k {W^*}^l$.  Now, if $f$
belongs to the $n$-spectral subspace, then $f=P_n(f)$, where $P_n$ is
the corresponding spectral projection.  On the other hand, if $P_n$ is
applied to the linear combination just mentioned, all terms will vanish
except for those for which $i-j+k-l=n$.

  The somewhat curious fact that $Z$ and $W$ are also eigenvalues for
$\t$ implies that
  $$
  Z^i \defprod {Z^*}^j \defprod W^k \defprod {W^*}^l =
  \mu Z^i {Z^*}^j W^k {W^*}^l,
  $$
 for some complex number $\mu$ of modulus one.  Therefore one concludes
that each spectral subspace for the deformed gauge action is contained
in the sub-algebra of
  $\DA{C(S^3)}{\g}{\t}$
  generated by $Z$ and $W$.  This shows that $Z$ and $W$ generate
  $\DA{C(S^3)}{\g}{\t}$
  and hence that $\psi$ is surjective.

  We next show that $\psi$ is injective.  Consider the circle action on
  $S^3_\tt$
  specified, on the generators, by
  $$
  \alpha_\lambda(B) = \lambda B \and
  \alpha_\lambda(C) = \lambda C,
  $$
  for $\lambda\in S^1$.  The homomorphism $\psi$, under scrutiny, is
clearly equivariant for the action just defined on
  $S^3_\tt$
  and the deformed gauge action $\tilde\g$ on
  $\DA{C(S^3)}{\g}{\t}$.
  By using
  \scite{\Circle}{2.9},
  it is now enough to verify that $\psi$ is injective on the fixed point
sub-algebra of $S^3_\tt$ for $\alpha$.  Let us denote that sub-algebra
by $F$.

  Recall that Matsumoto \scite{\MatsumotoHopf}{Theorem 6} has shown that
$F$ is isomorphic to the commutative \cstar-algebra of functions on the
two-sphere $S^2$.  More precisely, $F$ turns out to be generated by the
elements $M$ and $H$ of
  $S^3_\tt$
  given by
  $H = C^*C$
  and
  $M = CB^*$.  It is easy to see that these operators satisfy the
relations
  \zitemno = 0
  \zitem $H^* = H$,
  \zitem $M^*M = MM^*$, \hfill (\seqnumbering) \label{\SphereRelations}
  \zitem $MH = HM$, and
  \zitem $M^*M + H^2 = H$.

  Matsumoto has, in fact, shown that $F$ is the universal \cstar-algebra
on generators $H$ and $M$ satisfying the above relations.

  Now, the images of $H$ and $M$ under $\psi$ are, of course,
  $$
  \psi(H) = W\defstar\defprod W = W^*W
  $$
  and
  $$
  \psi(M) = W \defprod Z\defstar = W Z^*,
  $$
  both of which lie in the fixed point sub-algebra, say $B_0$, for the
deformed gauge action on
  $\DA{C(S^3)}{\g}{\t}$.  The crucial point is that this algebra is
impervious to the deformation, as one can easily see, so $B_0$ is just
the algebra of continuous functions on the quotient space $S^3/S^1$,
which is homeomorphic to $S^2$.

  An explicit homeomorphism between
  $S^3/S^1$ and $S^2$
  may be given by mapping the quotient class of $(z,w) \in S^3$ to
  the pair $(h,m) \in \R\x\C$, defined by
  $(h,m) = (w\bar w,w\bar z)$.  It is elementary to check that $(h,m)$
satisfies the equation
  $$
  |m|^2 + h^2 = h,
  $$
  which is precisely the equation defining the sphere of radius ${1\over
2}$ centered at $({1\over 2}, 0+i0)$ in $\R\x\C$.  The map
  $$
  \class{(z,w)} \vaipara (w\bar w,w\bar z)
  $$
  can now be shown to provide a homeomorphism from
  $S^3/S^1$
  onto the above mentioned model for the 2-sphere.

  Whenever a compact subset $K$ of $\R\x\C$ is defined via an equation
(or even a system of equations), such as the sphere above, one can prove
that $C(K)$ is the universal \cstar-algebra generated by symbols $h$ and
$m$, subject to the conditions
  \zitemno=0
  \zitem $h^* = h$,
  \zitem $m^*m = mm^*$,
  \zitem $mh = hm$,

  \medskip\noindent to which one should add the equations used to define
$K$.  This implies that $B_0$ is the universal \cstar-algebra generated
by a pair of elements (namely $h=W^*W$ and $m=WZ^*$) subject to the same
relations as the ones defining $F$, that is \lcite{\SphereRelations}.

  Therefore one sees that $\psi$ is an isomorphism between $F$ and
$B_0$, hence injective.  By \scite{\Circle}{2.9}, it follows that $\psi$
is injective everywhere and thus it is an isomorphism.
  \proofend

  In order to discuss the infinitesimal aspects of the deformation of
$S^3$ recently described, let $D_1$ and $D_2$ denote the differential
operators defined by
  $$
  D_1 f(z,w) = \deriv{f(\expi{\lambda}z,w)}{\lambda}{0},
  $$
  and
  $$
  D_2 f(z,w) = \deriv{f(z,\expi{\lambda}w)}{\lambda}{0},
  $$
  for $(z,w) \in S^3$ and $f\in C^\infty(S^3)$.

  \state Theorem \label{\DerivShpere} If $f$ and $g$ are in
  $C^\infty(S^3)$
  then
  $$
  \lim_{\tt \rightarrow 0}
  \left\Vert {f \defprod_\tt g - g \defprod_\tt f \over \tt}
  -
  {1 \over 2\pi i}
  \(D_2(f) D_1(g) - D_1(f) D_2(g)\)
  \right\Vert_\tt
  = 0,
  $$
  where $\defprod_\tt$ and $\|\cdot\|_\tt$ refer to the deformed
multiplication and norm of $S^3_\tt$.
  Therefore, the family
  $\{
  S^3_\tt
  \}_{\tt\in \R}$
  gives a strict deformation quantization for $C(S^3)$, in the direction
of the Poisson bracket $D_2\wedge D_1$.

  \proof Let $\action$ be the action of $\R^2$ on $S^3$ defined by
  $$
  \action_{(x,y)}(z,w) = (\expi{(x+y)}z, \expi{x}w)
  \for (x,y)\in \R^2,\ (z,w)\in S^3.
  $$

  As in section \lcite{\DerivSection} we may use $\action$ to obtain the
\dd $(\Z,\gamma,\t^\hbar)$.  However, one can easily see that this is
precisely the \dd used earlier in this section for $\hbar=\tt$.  So, we
may use \lcite{\PoissonCorollary} to treat the deformation $S^3_\tt$.
But, before that, let us remark that, since $\phi$ is a smooth action of
$\R^2$ on the compact manifold $S^3$, then any smooth function on $S^3$
will be $\action$-smooth.  Also, for $f$ in $C^\infty(S^3)$ we have,
using the notation of section \lcite{\DerivSection},
  $$
  \pdop{x}{f}\calcat{(z,w)} =
  \deriv{f(\expi{\lambda}z, \expi{\lambda}w)}{\lambda}{0}
  = D_1f(z,w) + D_2f(z,w),
  $$
  and
  $$
  \pdop{y}{f}\calcat{(z,w)} =
  \deriv{f(\expi{\lambda}z,w)}{\lambda}{0}
  = D_1f(z,w),
  $$
  that is,
  $\partial_x = D_1 + D_2$ while $\partial_y = D_1$. The Poisson bracket
appearing in \lcite{\PoissonCorollary} then becomes
  $$
  \pdop{x}{f} \pdop{y}{g} - \pdop{y}{f} \pdop{x}{g}
  =
  (D_1(f) + D_2(f)) D_1(g) -
  D_1(f)(D_1(g) + D_2(g))
  \$=
  D_2(f) D_1(g) - D_1(f) D_2(g),
  $$
  concluding the proof.
  \proofend

  \section {Example: Non commutative Lens spaces}
  Matsumoto and Tomiyama \cite{\MaTo}, building on \cite{\Mt}, have
introduced non-commutative versions of the classical lens spaces.  This
section is dedicated to proving that these can be described by using our
method of deformation.

  Recall that for nonzero co-prime integers $p$ and $q$, with $p\neq 0$,
the lens space $L(p,q)$ can be defined to be the quotient of the
three-sphere $S^3$ by the action of the finite cyclic group $\Z_p$
generated by the diffeomorphism
  \def\eigvalZ{\rho}
  \def\eigvalW{\rho^q}
  $$
  \tau(z,w) =
  (\eigvalZ z,\eigvalW w)
  \for (z,w)\in S^3,
  $$
  where $\rho =\expi{/p}$.

  Observe that, if one induces $\tau$ to an automorphism of $C(S^3)$ by
the formula
  $$
  \tau(f)\calcat{(z,w)} = f
  (\eigvalZ z,\eigvalW w)
   \for f\in C(S^3), \quad (z,w)\in S^3,
  $$
  then we have $\tau(Z) = \eigvalZ Z$ and $\tau(W) = \eigvalW W$, where
$Z $ and $W$ are the coordinate functions on $S^3$ (defined in the proof
of \lcite{\NCSphere}).
  Since $Z$ and $W$ generate $C(S^3)$, these equations actually define
$\tau$.  In addition one sees that the fixed point sub-algebra of
$C(S^3)$ for $\tau$ coincides with the algebra of continuous functions
on the quotient $S^3/\Z_p = L(p,q)$.

Let $\tt$ be a real number, fixed throughout.  Consider the automorphism
$\sigma$ of
  $S^3_\tt$
  prescribed by
  $\sigma(B) = \eigvalZ B$ and
  $\sigma(C) = \eigvalW C$,
  where $B$ and $C$ are as in the previous section.
  Among the many different characterizations of the non-commutative lens
space $L_\tt(p,q)$ presented in \cite{\MaTo}, one has:

  \definition The non-commutative lens space $L_\tt(p,q)$ is defined to
be the fixed point sub-algebra of $S^3_\tt$ under the automorphism
$\sigma$.

  Regarding the \dd $(\Z,\g,\t)$ for the algebra $C(S^3)$, defined
shortly before \lcite{\NCSphere}, in terms of a given value for the
parameter $\tt$, observe that $\g$ and $\t$ commute with $\tau$ and
hence both $\g$ and $\t$ leave invariant the fixed point sub-algebra for
$\tau$, which we have seen to be a model for $C(L(p,q))$.  By abuse of
language we still denote by $\g$ and $\t$ the corresponding restrictions
of these to $C(L(p,q))$.  So, this gives a \dd for $C(L(p,q))$ and we
may then form the deformed algebra $\DA{C(L(p,q))}{\g}{\t}$.

  \state Theorem For each real $\tt$ and co-prime integers $p$ and $q$,
with $p\neq 0$, the \cstar-algebras $\DA{C(L(p,q))}{\g}{\t}$ and
$L_\tt(p,q)$ are isomorphic.

  \proof We shall derive this from \lcite{\FixedPoint}.  In fact, let
the algebra $B$, mentioned in the statement of \lcite{\FixedPoint} be
$C(S^3)$ with the \dd $(\Z,\g,\t)$ referred to above.  Then, as we have
seen in \lcite{\NCSphere}, $\DA{C(S^3)}{\g}{\t}$ is naturally isomorphic
to $S^3_\tt$.  Still referring to the statement of \lcite{\FixedPoint},
let $H = \Z_p$, which acts on $B$ via $\tau$.  The natural extension of
$\tau$ to $\DA{C(S^3)}{\g}{\t}$, provided by \lcite{\TildeAction},
coincides with $\tau$ on the algebraic direct sum of the spectral
subspaces for the deformed gauge action and hence satisfies
  $$
  \tilde\tau(Z) = \eigvalZ Z
  \and
  \tilde\tau(W) = \eigvalW W.
  $$

  But, since the isomorphism between $\DA{C(S^3)}{\g}{\t}$ and $S^3_\tt$
puts $Z$ and $W$ in correspondence with $B$ and $C$, respectively, we
see that $\tilde\tau$ and $\sigma$ correspond to each other under this
isomorphism.  In particular, the fixed point sub-algebra for $\sigma$,
a.k.a $L_\tt(p,q)$, is isomorphic to the fixed point sub-algebra of
$\DA{C(S^3)}{\g}{\t}$ under $\tilde\tau$, which, by \lcite{\FixedPoint},
is isomorphic to the deformed algebra $\DA{C(L(p,q))}{\g}{\t}$.
  \proofend

  Observe that $\tau$ commutes with the action $\action$ referred to in
the proof of \lcite{\DerivShpere}.  Therefore the operators $D_1$ and
$D_2$ leave $C(L(p,q))$ invariant, when this is viewed as a subset of
$C(S^3)$.

  \state Theorem If $f$ and $g$ are in
  $C^\infty(L(p,q))$
  then
  $$
  \lim_{\tt \rightarrow 0}
  \left\Vert {f \defprod_\tt g - g \defprod_\tt f \over \tt}
  -
  {1 \over 2\pi i}
  \(D_2(f) D_1(g) - D_1(f) D_2(g)\)
  \right\Vert_\tt
  = 0,
  $$
  where $\defprod_\tt$ and $\|\cdot\|_\tt$ refer to the deformed
multiplication and norm of $L_\tt(p,q)$.
  Therefore, the family
  $\{
  L_\tt(p,q)
  \}_{\tt\in \R}$
  gives a strict deformation quantization for $C(L(p,q))$, in the
direction of the Poisson bracket $D_2\wedge D_1$.

  \proof This is, in view of the comment above, a direct application of
\lcite{\DerivShpere} for $f$ and $g$ in $C(L(p,q))$.
  \proofend

  \section {Example: Non commutative Heisenberg manifolds}
  For each positive integer $c$, the Heisenberg manifold $M^c$ consists
of the quotient $H/G^c$, where $H$ is the Heisenberg group
  $$
  H = \left\{ \pmatrix{1 & y & z \cr
                       0 & 1 & x \cr
                       0 & 0 & 1 \cr}
              : x,y,z \in \R \right\},
  $$
  viewed as a subgroup of $SL_3(\R)$,
  and $G^c$ is the discrete subgroup obtained when $x$, $y$ and $cz$ are
required to be whole numbers.

  To facilitate our notation we will use the coordinate system on $H$
suggested by the association
  $$
  (x,y,z)
  \leftrightarrow
  \pmatrix{1 & y & z \cr
  0 & 1 & x \cr
  0 & 0 & 1 \cr}
  $$
  and hence we will identify $H$ with $\R^3$ without further warning.
  Under this coordinate system the multiplication in $H$ becomes
  $$
  (x,y,z) (m,n,p) = (x+m,y+n,z+p+ym).
  \eqno{(\seqnumbering)}
  \label{\GroupLaw}
  $$
  So, $M^c$ can be described as the quotient of the Euclidean space
 $\R^3$ by the right action of $G^c$ given by \lcite{\GroupLaw}.

  In \cite{\RieffelHeis} Rieffel introduced a continuous field of
\cstar-algebras, denoted $D^c_{\mu,\nu}$, where $\mu$ and $\nu$ are real
parameters, such that
  $D^c_{\mu,\nu}$
  is isomorphic to
  $C(M^c)$
  when $\mu=\nu=0$.

  Recall from \cite{\AEE} that
  $D^c_{\mu,\nu}$
  can be defined to be the crossed product of $C(\T^2)$, the algebra of
continuous functions on the two-torus, by a certain Hilbert bimodule.
There are, in fact, many descriptions for this construction.  Perhaps
the simpler such is provided by
  \cite{\Picard}, where it is shown that
  $$
  D^c_{\mu,\nu} \simeq
  C(\T^2) \crossproduct_{X^c_{\alpha_{\mu,\nu}}} \Z.
  $$
  To describe this in detail, let
  $X^c$ be the set of
  continuous complex-valued functions on two real variables $x$ and $y$
satisfying
  \zitemno = 0
  \zitem $f(x,y+1) = f(x,y)$, and
  \zitem $f(x+1,y) = \expineg{cy}f(x,y)$.

  \medskip Viewing the elements of $C(\T^2)$ as periodic functions on
two real variables, it is easy to check that, under pointwise
multiplication,
  $C(\T^2) X^c \subseteq X^c$.
  In this way, $X^c$ is given the structure of a $C(\T^2)$--module.  If
we let, for $f$ and $g$ in $X^c$,
  $
  \<f,g\>_L = f\overline{g},
  $
  then it is clear that $\<f,g\>_L$ is a periodic function on both
variables, and hence belongs to $C(\T^2)$.  This gives $X^c$ the
structure of a left Hilbert module over $C(\T^2)$.  Defining
  $
  \<f,g\>_R = \overline{f}g,
  $
  $X^c$ becomes a Hilbert bimodule.

  Now, given real parameters $\mu$ and $\nu$, consider the automorphism
  $\alpha_{\mu,\nu}$
  of $C(\T^2)$
  given by
  $$
  \alpha_{\mu,\nu}(f)\calcat{(x,y)} = f(x+2\mu,y+2\nu).
  $$
  The Hilbert bimodule $X^c_{\alpha_{\mu,\nu}}$, appearing above, is
obtained by altering the right module operations of $X^c$ using
${\alpha_{\mu,\nu}}$ as in \cite{\Picard}, that is: if $X$ is any
Hilbert bimodule over a \cstar-algebra $A$ and $\alpha$ is an
automorphism of $A$, we let $X_\alpha$ denote the Hilbert bimodule over
$A$ which coincides with $X$ as a left Hilbert module but which is
equipped with the right module structure given by
  $$
  x \cdot a = x \alpha(a) \for x \in X,\ a \in A,
  $$
  and right inner product $\<\cdot,\cdot\>_R^{M_\alpha}$ given by
  $$
  \<x,y\>_R^{M_\alpha} = \alpha\inv(\<x,y\>_R) \for x,y \in X.
  $$

  So, $X^c_{\alpha_{\mu,\nu}}$ is a Hilbert bimodule over $C(\T^2)$ and
we may construct, as in \cite\AEE, the crossed product
  $C(\T^2) \crossproduct_{X^c_{\alpha_{\mu,\nu}}} \Z$ which, according
to \scite{\Picard}{Section 2}, is isomorphic to $D^c_{\mu,\nu}$.

  As in our earlier examples, we will show that
  $D^c_{\mu,\nu}$ can be described as a deformation of
  $C(M^c)$ relative to a certain \dd.

  Given $(x,y,z)$ in the Heisenberg group, we shall denote its class in
$H/G^c$ by $[x,y,z]$.

  Let $\mu$ and $\nu$ be fixed real numbers and consider the map
$\action$, from $\R^2$ into the Heisenberg group $H$, given by
  $$
  \action(a,b) \=
  \exp\(
  a\pmatrix{ 0 & 0 & 1/c \cr
             0 & 0 & 0 \cr
             0 & 0 & 0}
  +
  b\pmatrix{ 0 & 2\nu & 0 \cr
             0 & 0 & 2\mu \cr
             0 & 0 & 0}
  \)
  =
  \pmatrix{ 1 & 2b\nu & 2b^2\mu\nu+a/c \cr
            0 & 1 & 2b\mu \cr
            0 & 0 & 1}.
  $$
  Since the two summands being exponentiated commute with each other,
one sees that $\action$ is a group homomorphism, yielding an action of
$\R^2$ on $H$, by left multiplication.
  This action,
  which obviously commutes with the right action of $G^c$ on $H$,
  drops to the quotient, producing the following action of $\R^2$ on
$M^c$:
  $$
  \phi_{(a,b)}[x,y,z] =
  [x + 2b\mu, y + 2b\nu, z + 2b\nu x + 2b^2\mu\nu + a/c ],
  \eqno{(\seqnumbering)}
  \label{\GoodAction}
  $$
  for $(a,b)\in\R^2$ and $[x,y,z]\in M^c$.

  If we now let $(\Z,\g,\t^\hbar)$ be the \dd given by $\phi$, as in
section \lcite{\DerivSection}, the gauge action $\g$, seen, as before,
as an action of the circle group, will be
  $$
  \g_{\expi{t}}([x,y,z] ) = [x,y,z+t/c],
  $$
  while the deforming actions $\t^\hbar$ of
  $\Z$ on $M^c$ are given by iterating the diffeomorphism
  $$
  \t^\hbar([x,y,z]) = [x + 2\hbar\mu, y + 2\hbar\nu, z + 2\hbar\nu x +
2\hbar^2\mu\nu],
  \for [x,y,z]\in M^c.
  $$

  For the time being, let us assume that $\hbar=1$ or, what amounts to
the same, that $\mu$ and $\nu$ are replaced, respectively, by $\hbar\mu$
and $\hbar\nu$.  Correspondingly, let us denote the $\t^\hbar$ above
simply by $\t$.

  For each integer $k$ let us indicate by $B_k$ the $k$-spectral
subspaces for the gauge action $\g$ on $C(M^c)$. In particular, $B_0$,
which is the algebra of fixed points, coincides with the algebra of
continuous functions on the quotient $M^c/S^1$.  It is a simple task to
verify that the map
  $$
  [x,y,z] \in M^c \vaipara (\expi x,\expi y) \in \T^2
  $$
  drops to a homeomorphism from $M^c/S^1$ to the 2-torus $\T^2$.  In
other words, $B_0$ is isomorphic to $C(\T^2)$.

  In general, for each $k$ in $\Z$, the $k$-spectral subspace $B_k$ is
given by the set of functions $f : M^c \arw \C$ satisfying
$\g_\lambda(f) = \lambda^k f$, for $\lambda$ in $S^1$ or, equivalently,
  $$
  f[x,y,z+t/c] = \expi{kt}f[x,y,z].
  $$

  This reflects the fact that $\gamma$ is nothing but the dual action of
 $C(M^c)$, when the latter is viewed as a Hilbert--bimodule crossed
product \cite{\AEE}.
  Now, this implies that
  $
  f[x,y,z] = \expi{kcz}f[x,y,0].
  $
  So, $f$ is determined by its values on the elements $[x,y,0]$.  This
suggests defining, for each such $f$, the function
  $
  g(x,y) := f[x,y,0].
  $
  Since
  $(x,y,0)(0,1,0)=(x,y+1,0)$,
  we see that $g$ is periodic in its second variable.  Moreover, since
  $(x,y,0)(1,0,0)=(x+1,y,y)$,
  we have that
  $$
  g(x,y) = f[x,y,0] = f[x+1,y,y] = \expi{kcy} f[x+1,y,0] =
  \expi{kcy} g(x+1,y).
  $$
  Summarizing, we have
  \zitemno = 0
  \zitem $g(x,y+1) = g(x,y)$, and
  \zitem $g(x+1,y) = \expineg{kcy}g(x,y)$,

  \medskip\noindent
  which the reader should compare with the equations defining the
Hilbert bimodule $X^c$, earlier in this section.  Conversely, given any
continuous function $g : \R^2 \arw \C$ satisfying (i) and (ii) above,
one may define $f[x,y,z] = \expi{kcz}g(x,y)$, and, after verifying that
$f$ is indeed well defined, show that $f\in B_k$.

  We next observe that the gauge action on $M^c$ is semi-saturated,
  that is, $C(M^c)$ is generated, as a \cstar-algebra, by $B_0$ and
$B_1$.  This follows by the fact that $C(M^c)$ is a Hilbert--bimodule
crossed product \scite{\AEE}{Theorem 3.1}
 (see also \scite{\Circle}{4.1, 4.8} and \scite{\Amena}{6.2}).

  It is not hard to show, using the Tietze extension Theorem, that the
Fell bundle arising from any \stress{free} action of $\T^d$ on a locally
compact space, such as the one we have, is actually saturated.  However,
we will not need this fact presently.

  \state Theorem \label{\HeisIsom} For every integer $c$, and real
numbers $\mu$ and $\mu$, we have that $\DA{C(M^c)}{\g}{\t}$ is
isomorphic to $D^c_{\mu,\nu}$.

  \proof By \lcite{\Satura} we have that $\DB{\B}{\t}$ is also
semi-saturated and hence, by
 \cite{\AEE}, in conjunction with \scite{\Amena}{4.2 and 4.7}, we
conclude that
  $\DA{C(M^c)}{\g}{\t}$,
  that is
  $C^*(\DB{\B}{\t})$,
  is given by the Hilbert bimodule crossed product
  $
  B_0 \crossproduct_{B_1} \Z.
  $
  It is important to stress that
  the $B_0$--Hilbert bimodule structure of $B_1$ we are referring to, is
that coming from the operations of
  $C^*(\DB{\B}{\t})$, that is, the deformed bundle operations
  $\defprod$ and $\defstar$ of $\DB{\B}{\t}$.  To make this more
explicit, let $a\in B_0$ and $b,c\in B_1$, which we may assume are given
by
  $a[x,y,z] = f(x,y)$,
  $b[x,y,z] = \expi{cz}g(x,y)$, and
  $c[x,y,z] = \expi{cz}h(x,y)$, where $f$ is periodic and both $g$ and
$h$ satisfy the conditions (i) and (ii) above for $k=1$.  The reader may
then verify that
  $$
  \matrix{
  \hfill a\defprod b[x,y,z] &
         = & \expi{cz} f(x,y) g(x,y),\hfill\cr\cr
  \hfill b\defprod a[x,y,z] &
         = & \expi{cz} g(x,y) f(x+2\mu,y+2\nu),\hfill\cr\cr
  \hfill b\defstar\defprod c[x,y,z] &
         = & \overline{g(x-2\mu,y-2\nu)} h(x-2\mu,y-2\nu),\hfill\cr\cr
  \hfill b\defprod c\defstar[x,y,z] &
         = & g(x,y)\overline{h(x,y)}.\hfill\cr}
  $$
  These formulas tell us that the pair
  $(B_0,B_1)$
  is identical to
  $(C(\T^2), X^c_{\alpha_{\mu,\nu}})$
  as far as the Hilbert bimodule structure is concerned.  Hence, as we
already know that
  $\DA{C(M^c)}{\g}{\t}
      = B_0 \crossproduct_{B_1} \Z$
  and that
  $D^c_{\mu,\nu}
      = C(\T^2) \crossproduct_{X^c_{\alpha_{\mu,\nu}}} \Z$,
  it follows that
  $\DA{C(M^c)}{\g}{\t} \simeq D^c_{\mu,\nu}$.
  \proofend

  Let us now compute the differential operators $\partial_x$ and
$\partial_y$, as in section \lcite{\DerivSection}, arising from the
action $\phi$ of $\R^2$ on $M^c$.
  However, to avoid a notational conflict with the already established
coordinate system $[x,y,z]$ for $M^c$, we will denote them by
$\partial_a$ and $\partial_b$, respectively.  For a smooth function $f$
on $M^c$, we have
  $$
  \pdop{a}{f}[x,y,z] =
  \deriv{f[x,y,z+a/c]}{a}{0} =
  c\inv\pdop{3}{f}[x,y,z],
  $$
  while
  $$
  \pdop{b}{f}[x,y,z] =
  \deriv{f[x+2b\mu,y+2b\nu,z+2b^2\mu\nu+2b\nu x]}{b}{0}
  \$=
  \bigl(
  2\mu\pdop{1}{f} + 2\nu\pdop{2}{f} + 2\nu x\pdop{3}{f}
  \bigr)[x,y,z]
  $$
  where $\partial_1$, $\partial_2$ and $\partial_3$ correspond to the
partial differentiation operators for the standard coordinates on
$\R^3$.

The relevant Poisson bracket on $M^c$ becomes
  $$
  \{\cdot,\cdot\} =
  \partial_x\wedge\partial_y =
  c\inv \partial_3 \wedge
  (2\mu\partial_1 + 2\nu\partial_2 + 2\nu x\partial_3)
  =
  2c\inv \partial_3 \wedge
  (\mu\partial_1 + \nu\partial_2),
  $$
  which, up to a multiplicative factor, is precisely the Poisson bracket
of interest in section 2 of \cite{\RieffelHeis}.
  We may therefore deduce from \lcite{\PoissonCorollary} and
\lcite{\HeisIsom}, one of the main results of \cite{\RieffelHeis}:

  \state Theorem
  The family
  $\{
  D^c_{\hbar\mu,\hbar\nu}
  \}_{\hbar\in \R}$
  forms a strict deformation quantization in the direction of the
Poisson bracket
  $2c\inv \partial_3 \wedge(\mu\partial_1 + \nu\partial_2)$.

  \bigbreak
  \centerline{\tensc References}
  \nobreak\medskip\frenchspacing
  \catcode`\@=0

@TechReport{\AEE,
  auth = {B. Abadie, S. Eilers and R. Exel},
  title = {{M}orita Equivalence for Crossed Products by {H}ilbert
Bimodules},
  institution = {University of Copenhagen},
  year = {1994},
  type = {preprint},
  note = {to appear in {\sl Trans. Amer. Math. Soc\/}},
  NULL = {},
  akey = {AbadieEilersExel},
  author = {Beatriz Abadie, Soren Eilers and Ruy Exel},
  atrib = {IRA},
  }

@Unpublished{\Picard,
  auth = {B. Abadie and R. Exel},
  title = {Hilbert $C^*$-bimodules over commutative $C^*$-algebras and
an isomorphism condition for quantum {H}eisenberg manifolds},
  institution = {Universidade de S\~ao Paulo},
  year = {1996},
  type = {preprint},
  note = {to appear in {\sl Rev. Math. Phys\/}},
  NULL = {},
  akey = {AbadieExel},
  author = {Beatriz Abadie and Ruy Exel},
  atrib = {IRA},
  }

@Book{\Brat,
  auth = {O. Bratteli},
  title = {Derivations, Dissipations and Group Actions on C*-algebras},
  publisher = {Springer--Verlag},
  year = {1986},
  volume = {1229},
  series = {Lecture Notes in Mathematics},
  NULL = {},
  author = {Ola Bratteli},
  }

@Book{\Dixmier,
  auth = {J. Dixmier},
  title = {$C^*$-Algebras},
  publisher = {North Holland},
  year = {1982},
  volume = {},
  series = {},
  NULL = {},
  author = {J. Dixmier},
  }

@Article{\Soft,
  auth = {R. Exel},
  title = {The Soft Torus and Applications to Almost Commuting
Matrices},
  journal = {Pacific J. Math.},
  year = {1993},
  volume = {160},
  pages = {207--217},
  NULL = {},
  akey = {Exel},
  author = {Ruy Exel},
  number = {2},
  MR = {94f:46091},
  amsclass = {46L80 (46L05)},
  atrib = {IR},
  }

@Article{\Circle,
  auth = {R. Exel},
  title = {Circle Actions on {$C^*$}-Algebras, Partial Automorphisms and
a Generalized {P}imsner--{V}oiculescu Exact Sequence},
  journal = {J. Funct. Analysis},
  year = {1994},
  volume = {122},
  pages = {361--401},
  NULL = {},
  akey = {Exel},
  author = {Ruy Exel},
  number = {2},
  MR = {95g:46122},
  amsclass = {46L55 (46L80 47B35)},
  atrib = {IR},
  }

@Unpublished{\Amena,
  auth = {R. Exel},
  title = {Amenability for {F}ell Bundles},
  institution = {Universidade de S\~ao Paulo},
  year = {1996},
  type = {preprint},
  note = {},
  NULL = {},
  akey = {Exel},
  author = {Ruy Exel},
  atrib = {S},
  }

@Book{\FD,
  auth = {J. M. G. Fell and R. S. Doran},
  title = {Representations of *-algebras, locally compact groups, and
Banach *-algebraic bundles},
  publisher = {Academic Press},
  year = {1988},
  volume = {125 and 126},
  series = {Pure and Applied Mathematics},
  NULL = {},
  author = {J. M. G. Fell and R. S. Doran},
  }

@Article{\QCan,
  auth = {P. E. T. Jorgensen, L. M. Schmitt and R. F. Werner},
  title = {$q$-canonical commutation relations and stability of the
Cuntz algebra},
  journal = {Pacific J. Math.},
  year = {1994},
  volume = {165},
  pages = {131--151},
  NULL = {},
  note = {},
  author = {},
  MR = {95g:46116},
  amsclass = {46L40 (46L35 46L60 81S05)},
  number = {1},
  }

@Article{\Mt,
  auth = {K. Matsumoto},
  title = {Noncommutative three-dimensional spheres},
  journal = {Japan. J. Math. (N.S.)},
  year = {1991},
  volume = {17},
  pages = {333--356},
  NULL = {},
  author = {Kengo Matsumoto},
  }

@Article{\MatsumotoHopf,
  auth = {K. Matsumoto},
  title = {Noncommutative three-dimensional spheres. II.  Noncommutative
Hopf fibering},
  journal = {Yokohama Math. J.},
  year = {1991},
  volume = {38},
  pages = {103--111},
  NULL = {},
  author = {Kengo Matsumoto},
  number = {2},
  MR = {93c:46129},
  amsclass = {46L87 (58B30)},
  atrib = {},
  }

@Article{\MaTo,
  auth = {K. Matsumoto and J. Tomiyama},
  title = {Noncommutative lens spaces},
  journal = {J. Math. Soc. Japan},
  year = {1992},
  volume = {44},
  pages = {13--41},
  NULL = {},
  author = {Kengo Matsumoto and Jun Tomiyama},
  }

@TechReport{\Ng,
  auth = {C.-K. Ng},
  title = {Reduced Cross-sectional $C^*$-algebras of $C^*$-algebraic
bundles and Coactions},
  institution = {Oxford University},
  year = {1996},
  type = {preprint},
  note = {},
  NULL = {},
  author = {Chi-Keung Ng},
  }

@Article{\Irrat,
  auth = {M. A. Rieffel},
  title = {$C^*$-algebras associated with irrational rotations},
  journal = {Pacific J. Math.},
  year = {1981},
  volume = {93},
  pages = {415--429},
  NULL = {},
  author = {Marc A. Rieffel},
  }

@Article{\RieffelCF,
  auth = {M. A. Rieffel},
  title = {Continuous fields of $C^*$-algebras coming from group
cocycles and actions},
  journal = {Math. Ann.},
  year = {1989},
  volume = {283},
  pages = {631-643},
  NULL = {},
  author = {Marc A. Rieffel},
  }

@Article{\RieffelHeis,
  auth = {M. A. Rieffel},
  title = {Deformation quantization of Heisenberg manifolds},
  journal = {Commun. Math. Phys.},
  year = {1989},
  volume = {122},
  pages = {531--562},
  NULL = {},
  author = {Marc A. Rieffel},
  }

@Article{\RieffelDQi,
  auth = {M. A. Rieffel},
  title = {Deformation quantization for actions of ${\bf R}\sp d$},
  journal = {Mem. Amer. Math. Soc.},
  year = {1993},
  volume = {106},
  pages = {93 pp},
  NULL = {},
  author = {Marc A. Rieffel},
  number = {506},
  MR = {94d:46072},
  amsclass = {46L87 (58B30 81R50 81S10)},
  }

@Article{\RieffelDQii,
  auth = {M. A. Rieffel},
  title = {Quantization and $C\sp *$-algebras},
  journal = {Contemp. Math.},
  year = {1994},
  volume = {167},
  pages = {66--97},
  NULL = {},
  note = {$C\sp *$-algebras: 1943--1993 (San Antonio, TX, 1993)},
  author = {Marc A. Rieffel},
  MR = {95h:46108},
  amsclass = {46L60 (81S10)},
  }

@Article{\Woro,
  auth = {S. L. Woronowicz},
  title = {Twisted $SU_2$ groups. An example of a non-commutative
differential calculus},
  journal = {Publ. RIMS, Kyoto Univ.},
  year = {1987},
  volume = {23},
  pages = {117--181},
  NULL = {},
  author = {S. L. Woronowicz},
  }

@Article{\ZM,
  auth = {G. Zeller-Meier},
  title = {Produits crois\'es d'une $C^*$-alg\`ebre par un group
d'automorphismes},
  journal = {J. Math Pures Appl.},
  year = {1968},
  volume = {47},
  pages = {101-239},
  NULL = {},
  author = {G. Zeller-Meier},
  }

  \catcode`\@=12
  
  \begingroup
  \parindent=0pt
  \bigskip
  \obeylines
  
  Centro de Matem\'aticas
  Facultad de Ciencias
  Eduardo Acevedo 1139
  CP 11200, Montevideo -- Uruguay
  abadie@cmat.edu.uy

  \bigskip

  Departamento de Matem\'atica
  Universidade de S\~ao Paulo
  Rua do Mat\~ao, 1010
  05508-900 S\~ao Paulo -- Brazil
  exel@ime.usp.br

  \endgroup
  \medskip
  \rightline{May 1997}

  \bye